\pdfoutput=1
\documentclass[runningheads]{llncs}
\usepackage{amssymb}
\usepackage{amsmath}
\usepackage{booktabs} 
\usepackage{graphicx}
\usepackage{pgfplots}
\usepackage[all]{nowidow}
\usepackage[utf8]{inputenc}
\usepackage{tikz}
\usetikzlibrary{er,positioning,bayesnet}
\usepackage{multicol}
\usepackage{algpseudocode,algorithm,algorithmicx}
\usepackage{multirow}
\usepackage[inline]{enumitem} 
\usepackage[a4paper,left=32mm,
 top=20mm,right=32mm]{geometry}

\usepackage{fourier} 
\usepackage{ upgreek } 

\usepackage{colortbl} 
 \usepackage{color}
\definecolor{Silver}{rgb}{0.752,0.752,0.752}
\usepackage{hyperref}

\definecolor{blue}{HTML}{1F77B4}
\definecolor{orange}{HTML}{FF7F0E}
\definecolor{green}{HTML}{2CA02C}

\pgfplotsset{compat=1.14}

\allowdisplaybreaks

\begin{document}
\title{Wasserstein Distributionally Robust Optimization with Expected Value Constraints}
\titlerunning{DRWO with Expected Constraints}
%
\author{Diego Fonseca \and
Mauricio Junca }
%
%
\institute{Department of Mathematics, Universidad de los Andes, Bogot\'a, Colombia
\email{\{df.fonseca,mj.junca20\}@uniandes.edu.co}}
\maketitle              
\begin{abstract}
We investigate a stochastic program with expected value constraints, addressing the problem in a general context through Distributionally Robust Optimization (DRO) approach using Wasserstein distances, where the ambiguity set depends on the decision variable. We demonstrate that this approach can be reformulated into a finite-dimensional optimization problem, which, in certain instances, can be convex. Moreover, we establish criteria for determining the feasibility of the problem concerning the Wasserstein radius and the parameter governing the constraint. Finally, we present numerical results within the context of portfolio optimization. In particular, we highlight the distinctions between our approach and several existing non-robust methods, using both simulated data and real financial market data.

\end{abstract}
\section{Introduction}

In this work, we examine stochastic programs with expected value constraints, given by the following formulation:
\begin{equation}\label{StochsticProgWithExpectConst}
J=\left\{\begin{array}{ll}
{\displaystyle \min_{x\in\mathbb{R}^{m}} } & \Phi\left(F(x,\xi),\xi\sim\mathbb{P}\right)\\
\mbox{subject to} & \mathbb{E}_{\xi\sim\mathbb{P}}\left[G\left( x,\xi\right)\right]\geq \mu,\\
& x\in\mathcal{X},
\end{array}  \right.
\end{equation}
where $F$ and $G$ are functions such that $F,G:\mathbb{R}^{m}\times \mathbb{R}^{n}\rightarrow\mathbb{R}$, $\xi\in\mathbb{R}^{n}$ is a random vector with an (unknown) probability distribution $\mathbb{P}$ supported in $\Xi\subseteq\mathbb{R}^{n}$, that is, $\mathbb{P}\in\mathcal{P}(\Xi)$, and $\mathcal{X}\subseteq \mathbb{R}^{m}$ represents a set of constraints on the decision vectors. Additionally, the objective function $\Phi$ is a risk function dependent on the performance function $F$. For example, $\Phi$ can be the variance or the expected value of $F(x,\xi)$ with respect to $\xi$. The objective is to propose a strategy that addresses this problem, assuming that $\mathbb{P}$ is unknown but realizations of the random vector $\xi$ are available.

In cases where $\Phi\left(F(x,\xi),\xi\sim\mathbb{P}\right):=\mathbb{E}_{\xi\sim\mathbb{P}}\left[F(x,\xi)\right]$, this problem emerges in various contexts such as finance \cite{Li2017,Dentcheva2003}, operations research \cite{Miller1965}, and machine learning \cite{Rigollet2011,Mu2017}. Many attempts to solve (\ref{StochsticProgWithExpectConst}) employ Sample Average Approximation (SAA) \cite{WangAhmed2008,ShapiroMonteCarlo2003}, utilizing samples of $ \xi$ to replace expected values with sample means. Alternative strategies include those based on stochastic gradient or subgradient descent methods, as applied in \cite{Akhtar2021,Lan2020,Xiao2019,ZhangZhang2022}. These approaches can be sensitive to changes in sample quality, and out-of-sample performance may not always be optimal, especially as constraints might not be satisfied out-of-sample with small sample sizes. This is due to the need for at least one sample point during each iteration. Consequently, in order to achieve satisfactory optimality, larger samples may be necessary. Lastly, it is worth mentioning that other risk functions, such as Conditional Value-at-Risk, can also be formulated as (\ref{StochsticProgWithExpectConst}), as demonstrated in \cite{Rockafellar2000}.

When $\Phi\left(F(x,\xi),\xi\sim\mathbb{P}\right):=\mathrm{Var}_{\xi\sim\mathbb{P}}\left[F(x,\xi)\right]$, (\ref{StochsticProgWithExpectConst}) is referred to as the mean-variance model. This problem has been predominantly explored in portfolio optimization and inventory management, with one of its earliest appearances in \cite{Markowitz1952}. Consequently, the analysis of this problem has focused on specific $F$ functions and therefore some strategies employed in those contexts may not be applicable to problems with general $F$ functions. Nevertheless, most strategies used to address (\ref{StochsticProgWithExpectConst}) are concentrated in the set of strategies employed for the same problem in portfolio optimization. Therefore, it is essential to review the strategies used in that context.

In that case, $m=n$ where $m$ is the number of portfolio assets, $\xi$ is a random vector of returns for each asset, $x$ is a portfolio weights vector, and other constraints admissible for the investor are described by the set $\mathcal{X}$. Furthermore, in this scenario, $F=G$ where $F(x,\xi):= \langle x,\xi\rangle$ with $\langle\cdot,\cdot\rangle$ representing the Euclidean inner product in $\mathbb{R}^{m}$, such that $\langle x,\xi\rangle$ is the return of the portfolio $x$. Initial attempts to solve this problem considered estimates of the vector of means and the covariance matrix of returns. However, \cite{Chopra1993} demonstrated that the resulting portfolios perform poorly out-of-sample and are highly sensitive to variations in the estimates.

To address this issue, one of the first ideas was to treat the vector of means and the covariance matrix as variables, meaning the optimization problem's variables include the portfolio weights, the vector of means, and the covariance matrix. The choice of the feasible set for the vector and the matrix is crucial in this approach. Some works utilize sets based on a priori information about the returns or impose sets that are computationally tractable, as demonstrated in \cite{El-Ghaoui2003}, \cite{Zymler2011}, \cite{Natarajan2010}, \cite{Lotf2017}, \cite{Lotf2018}, and \cite{Won2020}. However, imposing unverifiable assumptions about the moments of the returns can also impact the out-of-sample performance of these methods. Moreover, the use of the vector of means and the covariance matrix arises naturally in problem (\ref{StochsticProgWithExpectConst}) for this portfolio optimization case since these two terms appear directly in the objective function and the constraint. 

In light of the described scenario, a data-driven approach for addressing (\ref{StochsticProgWithExpectConst}) based on Distributionally Robust Optimization (DRO) emerges, which involves solving the following optimization problem:
\begin{equation}\label{MarkovizRobustGeneral}
\left\{\begin{array}{ll}
{\displaystyle \min_{x\in\mathbb{R}^{m}} } & {\displaystyle\sup_{\mathbb{Q}\in\mathcal{D} }  \Phi\left(F(x,\xi),\xi\sim\mathbb{Q}\right)} \\
\mbox{subject to} & {\displaystyle\inf_{\mathbb{Q}\in\mathcal{D} }  \mathbb{E}_{\xi\sim\mathbb{Q}}\left[G\left( x,\xi\right) \right]   \geq \mu},\\
&  x\in\mathcal{X}.
\end{array}  \right.
\end{equation}
where $\mathcal{D}$ is a set of probability distributions, known as the \textit{ambiguity set}. However, DRO approach was initially proposed for unconstrained stochastic problems such as
\begin{equation}\label{ProbleEstocastico}
J_{\mathcal{D}}:=\min_{x\in\mathcal{X}}\mathbb{E}_{\xi\sim\mathbb{P}}[f(x,\xi)].
\end{equation}
where $f(x,\xi)$ is a cost function. In this context, the DRO approach for problem (\ref{ProbleEstocastico}) is formulated as
\begin{equation}\label{Eqn:DROprinsipal}
J_{\mathcal{D}}:=\min_{x\in\mathcal{X}}\sup_{\mathbb{Q}\in\mathcal{D}}\mathbb{E}_{\xi\sim\mathbb{Q}}[f(x,\xi)].
\end{equation}
Nonetheless, the results obtained to address (\ref{Eqn:DROprinsipal}) serve as the foundation for addressing (\ref{MarkovizRobustGeneral}), making it crucial to have a clear understanding of (\ref{Eqn:DROprinsipal}).

The selection of set $\mathcal{D}$ is vital for the tractability of (\ref{MarkovizRobustGeneral}), so it is also essential for (\ref{Eqn:DROprinsipal}). Approaches for defining $\mathcal{D}$ encompass sets of single-point supported distributions, distributions with moment restrictions, or parametric distribution families. Another method involves defining $\mathcal{D}$ as a ball centered on an empirical distribution using a distance notion, ensuring the true distribution $\mathbb{P}$ belongs to the ball with high probability or leads to decisions with satisfactory out-of-sample performance. For further information on this topic, we refer the reader to \cite{fonsecaDecDep2023} and the cited references. In this study, we adopt the Wasserstein distance, defining $\mathcal{D}$ as a ball centered at an empirical distribution with an appropriately chosen radius.

\begin{definition}[Wassertein distance]\label{Def:MetricaWasserstein}
The \textit{Wasserstein distance} $W_{p}(\mu,\nu)$ between $\mu,\nu\in\mathcal{P}_{p}(\Xi)$ is defined by
\begin{equation*}
{\displaystyle W_{p}(\mu,\nu):=\left(\inf_{\Pi\in\mathcal{P}(\Xi\times\Xi)}\left\{\int_{\Xi\times\Xi}\mathbf{d}^{p}(\xi,\zeta)\Pi(d\xi,d\zeta)\: :\: \Pi(\cdot \times\Xi)=\mu(\cdot),\: \Pi(\Xi\times\cdot)=\nu(\cdot)\right\}\right)^{1/p}
}
\end{equation*}
where 
$$\mathcal{P}_{p}(\Xi):=\left\{\mu\in\mathcal{P}(\Xi)\: :\: \int_{\Xi}\mathbf{d}^{p}(\xi,\zeta_{0})\mu(d\xi) < \infty\ \mbox{for some }\zeta_{0}\in\Xi\right\}$$
and $d$ is a metric in  $\Xi$.
\end{definition}
$W_{p}$ defines a metric in $\mathcal{P}_{p}(\Xi)$ for  $p\in[1,\infty)$, hence, the  ball with respect to some  $p$-Wasserstein distance with radius $\varepsilon>0$ and center $\mu\in \mathcal{P}(\Xi)$ is given by
\begin{equation}\label{BolaRespectoP}
\mathcal{B}_{\varepsilon}\left(\mu\right):=\left\{\nu\in \mathcal{P}(\Xi) \: \left|\: W_{p}(\mu,\nu)\leq\varepsilon \right.\right\}.
\end{equation}
Using $p$-Wasserstein distances, one could consider $\mathcal{D}=\mathcal{B}_{\varepsilon}\left(\widehat{\mathbb{P}}_{N}\right)$ in (\ref{Eqn:DROprinsipal}), which leads to the following formulation:
\begin{equation} \label{DROWGeneral}
\widehat{J}^{\mathrm{S}}_{N,p,q}(\varepsilon):=\min_{x\in\mathbb{X}}\sup_{\mathcal{Q}\in\mathcal{B}_{\varepsilon}\left(\widehat{\mathbb{P}}_{N}\right) }\mathbb{E}_{\xi\sim\mathbb{Q}}[f(x,\xi)],
\end{equation}
where the cost function used is $\mathbf{d}=\left\|\cdot\right\|_{q}$ (see definition \ref{Def:MetricaWasserstein}). In this work, we employ the notation ``S'' to denote  the ``standard'' formulation for distributionally robust optimization problems using Wasserstein's distance. Finally, the problem (\ref{DROWGeneral}) can be reformulated by the following theorem.

\begin{theorem}\label{Thm:ReformulacionDROWInterno}
Assume that $ f $  is upper semicontinuous with respect to $\xi$. Then the problem (\ref{DROWGeneral}) is equivalent to the optimization problem
\begin{equation}\label{Eqn:ReformulacionDROW}
\left\{
\begin{array}{lll}
{\displaystyle \inf_{x\in\mathcal{X},\lambda,s}} & {\displaystyle \lambda \varepsilon^{p} +\frac{1}{N}\sum_{i=1}^{N}s_{i}} &\\
\mbox{subject to} & {\displaystyle \sup_{\xi\in\Xi}\left(f(x,\xi)-\lambda \mathbf{d}^{p}(\xi,\widehat{\xi}_{i}) \right) \leq s_{i}  } & \forall i=1,\ldots,N, \\
&\lambda \geq 0.& 
\end{array}
\right.
\end{equation}
\end{theorem}
This theorem is established and proved in \cite{Blanchet2019}. Nonetheless, the reformulation (\ref{Eqn:ReformulacionDROW}) has been achieved under more stringent conditions in \cite{Kuhn2018} and \cite{Mehrotra2017}.

\begin{figure}[t] 
 \centering
\begin{tabular}{cc}
    \includegraphics[scale=0.31]{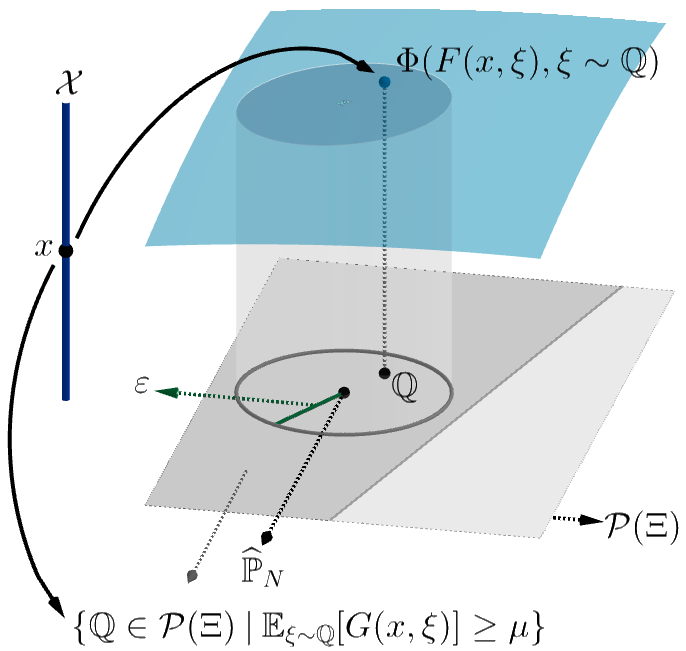} &  \includegraphics[scale=0.31]{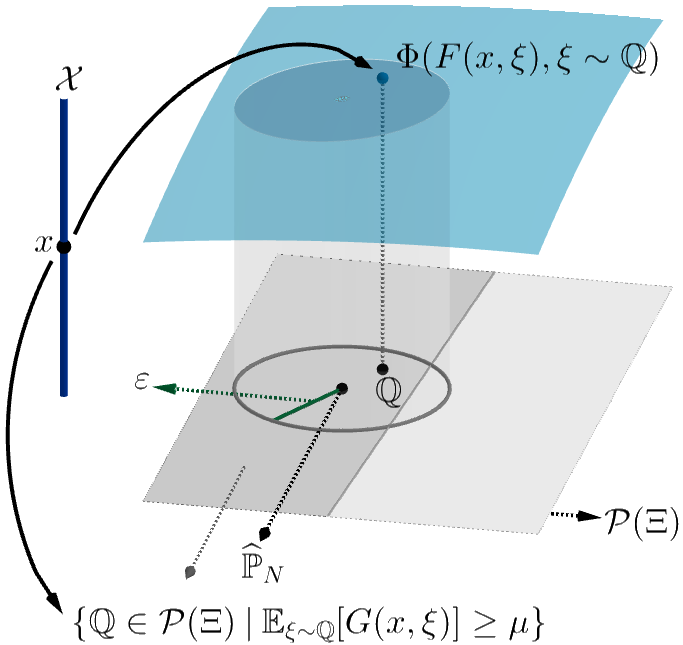}  \\
     (a) $x$ feasible. & (b) $x$ not feasible.
 \end{tabular} 
\caption{Illustration of the optimization process proposed in (\ref{MarkovizRobustGeneralWass}). } \label{fig:DiagramaStandard}
\end{figure}

The previous result suggests that it may be advantageous to consider $\mathcal{B}_{\varepsilon}\left(\widehat{\mathbb{P}}_{N}\right)$ as the ambiguity set in (\ref{MarkovizRobustGeneral}). Consequently, the following problem is obtained:
\begin{equation}\label{MarkovizRobustGeneralWass}
\widehat{J}_{N,p,q}^{\mathrm{S,cst}}(\varepsilon)=\left\{\begin{array}{ll}
{\displaystyle \min_{x\in\mathbb{R}^{m}} } & {\displaystyle\sup_{\mathbb{Q}\in\mathcal{B}_{\varepsilon}\left(\widehat{\mathbb{P}}_{N}\right) }  \Phi\left(F(x,\xi),\xi\sim\mathbb{Q}\right)} \\
\mbox{subject to} & {\displaystyle\inf_{\mathbb{Q}\in\mathcal{B}_{\varepsilon}\left(\widehat{\mathbb{P}}_{N}\right) }  \mathbb{E}_{\xi\sim\mathbb{Q}}\left[G\left( x,\xi\right) \right]   \geq \mu},\\
&  x\in\mathcal{X}.
\end{array}  \right.
\end{equation}
We will refer to this approach as the standard approach with constraints, denoted by the abbreviation ``S,cst''. Figure \ref{fig:DiagramaStandard} displays a diagram illustrating the optimization process representing (\ref{MarkovizRobustGeneralWass}), in which two cases can be observed: the first is when $x$ is feasible and the second is when it is not. This figure also shows that the decision variable $x$ influences the probability measure space $\mathcal{P}(\Xi)$ by inducing a feasible region and affecting the objective function.

Taking this into account, Theorem \ref{Thm:ReformulacionDROWInterno} is useful for acquiring a reformulation of (\ref{MarkovizRobustGeneralWass}).
For instance, if $\Phi\left(F(x,\xi),\xi\sim\mathbb{P}\right):=\mathbb{E}_{\xi\sim\mathbb{P}}\left[F(x,\xi)\right]$ and $F$ and $G$ are upper semicontinuous, (\ref{MarkovizRobustGeneral}) is equivalent to 
\begin{equation}\label{Eqn:ReformulacionDROWconConstraintRn}
\left\{
\begin{array}{lll}
{\displaystyle \inf_{x\in\mathcal{X},\lambda_{1},\lambda_{2},s,t,}}  &  {\displaystyle \frac{1}{N}\sum_{i=1}^{N}s_{i}}+ \lambda_{1} \varepsilon^{p} &\\
\mbox{subject to} & {\displaystyle \frac{1}{N}\sum_{i=1}^{N}t_{i}} -\lambda_{2} \varepsilon^{p} \geq \mu, & \\
&{\displaystyle \sup_{\xi\in\Xi}\left(F(x,\xi)-\lambda_{1} \mathbf{d}^{p}(\xi,\widehat{\xi}_{i}) \right) \leq s_{i}  } & \forall i=1,\ldots,N, \\
&{\displaystyle \inf_{\xi\in\Xi}\left(G(x,\xi)+\lambda_{2} \mathbf{d}^{p}(\xi,\widehat{\xi}_{i}) \right) \geq t_{i}  } & \forall i=1,\ldots,N, \\
&\lambda_{1},\lambda_{2} \geq 0.& 
\end{array}
\right.
\end{equation}
In some cases, depending on the forms of functions $F$ and $G$, this problem may present itself as a semi-infinite optimization problem with a large number of variables to optimize, which can be problematic. This occurs because the supremum or infimum appearing in the constraints of (\ref{Eqn:ReformulacionDROWconConstraintRn}) may not have explicit solutions. Solving semi-infinite programs is known to be computationally demanding. Additionally, there are several cases where the supremum and infimum in the constraints can be explicitly solved, but the resulting problem is either non-convex or computationally complex. This prompts us to propose an alternative type of ambiguity set.

The strategy involves using two distinct ambiguity sets, one in the objective function and another in the constraint. Both sets of ambiguity are defined according to (\ref{BolaRespectoP}); however, their radius and center depend on the decision variable. As a result, we consider decision-dependent ambiguity sets. This proposal will be explained in more detail in Section \ref{Sec:MeanVariance}. However, employing decision-dependent ambiguity sets is relatively uncommon in stochastic problems with stochastic constraints, as they predominantly appear in unconstrained stochastic problems. To appreciate the significance of using such ambiguity sets, readers may consult \cite{fonsecaDecDep2023} and the references included in that publication. In that sense, our proposal is based on \cite{fonsecaDecDep2023}.

Another important aspect to emphasize is the choice between addressing stochastic problems with constraints and dealing with regularized stochastic problems. Regularized problems typically do not involve constraints. To create a regularized version of (\ref{StochsticProgWithExpectConst}), the constraint function is incorporated into the objective function by multiplying it with a regularizing parameter and subsequently adding the result to the objective function. As a result, the regularizing parameter acts as a penalty, partially substituting the parameter $\mu$. Nevertheless, the value assigned to this parameter lacks the practical interpretability of the parameter $\mu$ in (\ref{StochsticProgWithExpectConst}). In the context of portfolio optimization, for example, $\mu$ is understood as the minimum return acceptable to the investor. In contrast, a regularizing parameter's value, depending on its magnitude, signifies the investor's level of preference or aversion to high returns. However, this preference level cannot be considered a general measure that can be replicated across different environments. If the data distribution changes, the initially chosen level may not necessarily produce decisions with the same performance as those generated before the distribution alteration. This issue is not present with $\mu$. However, it is important to acknowledge that regularized problems are generally computationally more manageable.


In summary, this work applies distributionally robust optimization (DRO) with the Wasserstein metric to address stochastic programs with expected value constraints, using a sample of the random vector. Our contributions can be outlined as follows:

\begin{enumerate}
\item[$\bullet$] We present a data-driven robust formulation of (\ref{StochsticProgWithExpectConst}) that emphasizes expected value constraints without the need for regularization parameters.
\item[$\bullet$] We demonstrate that our DRO approach for (\ref{StochsticProgWithExpectConst}) can be reformulated as optimization problems with finite-dimensional variables. In the context of portfolio optimization, we show that for specific cases, the problem resulting from our approach is convex.
\item[$\bullet$] We establish criteria to ensure the feasibility of the proposed approach, which involve the radius of the ball defining the ambiguity set, the level $\mu$ of the constraint, and the shape of the $F$ and $G$ functions. By identifying the sets of values for the radius and $\mu$ that make the proposed approach feasible, we can more effectively determine parameters that are likely to ensure good out-of-sample performance.
\item[$\bullet$] In portfolio optimization, we evaluate the performance of our approach and compare it with other traditional methods, using both synthetically generated return data and real market data. Our proposal demonstrates its advantages by achieving the highest expected return, cumulative wealth, and Sharpe ratio compared to other benchmarks, as well as a low turnover relative to the SAA strategy.
\end{enumerate}

The paper is structured as follows. In Section \ref{Sec:MeanVariance}, we describe our distributionally robust optimization model using the Wasserstein distance. We also derive tractable reformulations for the optimization problem and examine its feasibility. Additionally, we establish a criterion for estimating the probability that decisions generated by our approach satisfy the constraint of (\ref{StochsticProgWithExpectConst}) out-of-sample, which is expressed in terms of the ambiguity set's size. Simulation analysis of the proposed approaches is presented in Section \ref{Sec:PortfolioSubsection} for the portfolio optimization context. Finally, we draw conclusions in Section \ref{Sec:Conclusions}.

\sloppy \subsubsection*{Notation:} For  $q\in[1,\infty)\cap \mathbb{N}$, the $q$-norm in $\mathbb{R}^{k}$ is is noted as $\left\|\cdot\right\|_{q}$. For $N\in \mathbb{N}$, we let $[N] := \{1, 2,\ldots,N\}$. Additionally, the  $q$-Lipschitz norm of the function $f:\mathbb{R}^{n}\rightarrow \mathbb{R}$ is $\left\|f\right\|_{\mathrm{Lip,q}}:=\sup_{x\neq y}(f(x)-f(y))/\left\|x-y\right\|_{q}$. Finally, given a sample $\widehat{\xi}_{1},\ldots,\widehat{\xi}_{N}$ of a random variable $\xi$, the \textit{empirical distribution} of $\xi$ with respect to this sample is defined as the probability measure given by $\widehat{\mathbb{P}}_{N}:=\frac{1}{N}\sum_{i=1}^{N}\delta_{\widehat{\xi}_{i}}$, where $\delta_{x}$ is the Dirac delta function supported at $x$.


\section{Problem formulation and main results} \label{Sec:MeanVariance}

In this section, we introduce our approach to addressing (\ref{StochsticProgWithExpectConst}), which relies on defining the ambiguity set in (\ref{MarkovizRobustGeneral}) as proposed in \cite{fonsecaDecDep2023}. With this in mind, we must first consider the following assumption.

\begin{assumption}[Lipschitz]\label{AssumptionPrincipal}
We assume that $F$ and $G$ are $q$-Lipschitz functions with respect to $\xi$. This is, for each $x$, there exists $\gamma_{x,F,q} >0$ and $\gamma_{x,G,q} >0$ such that $|F(x,\xi)-F(x,\zeta)|\leq \gamma_{x,F,q}\left\|\xi-\zeta\right\|_{q} $ and $|G(x,\xi)-G(x,\zeta)|\leq \gamma_{x,G,q}\left\|\xi-\zeta\right\|_{q} $ for all $\xi,\zeta\in\Xi$ respectively. We denote $\gamma_{x,F,q}=\left\|F(x,\cdot)\right\|_{\mathrm{Lip,q}}$ and $\gamma_{x,G,q}=\left\|G(x,\cdot)\right\|_{\mathrm{Lip,q}}$.
\end{assumption}

Our approach incorporates an empirical distribution that depends on $x$. First, we introduce the following notation: For $x \in \mathbb{R}^{m}$, we define $\zeta^{x,F} := F(x, \xi)$ and $\zeta^{x,G} := G(x, \xi)$, noting that these are random variables. We denote the probability distributions of $\zeta^{x,F}$ and $\zeta^{x,G}$ as $\mathbb{P}^{x,F}$ and $\mathbb{P}^{x,G}$, respectively. Since they depend on $\mathbb{P}$, $\mathbb{P}^{x,F}$ and $\mathbb{P}^{x,G}$ are also unknown. Additionally, we define $\widehat{\zeta}^{x,F}_{i} := F(x, \widehat{\xi}_{i})$ and $\widehat{\zeta}^{x,G}_{i} := G(x, \widehat{\xi}_{i})$, such that $\widehat{\zeta}^{x,F}_{1}, \ldots, \widehat{\zeta}^{x,F}_{N}$ form a sample of $\zeta^{x,F}$, and $\widehat{\zeta}^{x,G}_{1}, \ldots, \widehat{\zeta}^{x,G}_{N}$ form a sample of $\zeta^{x,G}$. This enables us to define the empirical distributions of $\zeta^{x, F}$ and $\zeta^{x, G}$, given by $\widehat{\mathbb{P}}^{x,F}_{N} := \frac{1}{N} \sum_{i=1}^{N} \delta_{\widehat{\zeta}^{x,F}_{i}}$ and $\widehat{\mathbb{P}}^{x,G}_{N} := \frac{1}{N} \sum_{i=1}^{N} \delta_{\widehat{\zeta}^{x,G}_{i}}$. This dependence on $x$ is justified by the fact that the decision vector $x$ influences whether the constraint of (\ref{StochsticProgWithExpectConst}) is satisfied. Specifically, the constraint $\mathbb{E}_{\mathbb{P}}\left[G\left( x, \xi\right)\right] \geq \mu$ must be satisfied by $G\left( x, \xi\right)$, which depends on $x$, making it natural to introduce the dependence of $x$ in the ambiguity set. Consequently, we consider the following optimization problem for a given $\epsilon > 0$:

\begin{equation}\label{MarkovizRobust}
\widehat{J}_{N,p,q}^{\mathrm{A,cst}}(\varepsilon):=\left\{\begin{array}{ll}
{\displaystyle \min_{x\in\mathbb{R}^{m}} } & {\displaystyle\sup_{\mathbb{Q}\in\mathcal{B}_{\varepsilon\gamma_{x,F,q}}\left(\widehat{\mathbb{P}}_{N}^{x,F}\right) } \Phi\left(\zeta,\zeta\sim\mathbb{Q}\right) } \\
\mbox{subject to} & {\displaystyle\inf_{\mathbb{Q}\in\mathcal{B}_{\varepsilon\gamma_{x,G,q}}\left(\widehat{\mathbb{P}}_{N}^{x,G}\right) }  \mathbb{E}_{\mathbb{Q}}\left[\zeta\right]   \geq \mu},\\
&  x\in\mathcal{X}.
\end{array}  \right.
\end{equation}
where $\mathcal{B}_{\varepsilon\gamma_{x,F,q}}\left(\widehat{\mathbb{P}}_{N}^{x,F}\right)$ and $\mathcal{B}_{\varepsilon\gamma_{x,G,q}}\left(\widehat{\mathbb{P}}_{N}^{x,G}\right)$ represent balls centered at $\widehat{\mathbb{P}}_{N}^{x,F}$ and $\widehat{\mathbb{P}}_{N}^{x,G}$ with radii $\varepsilon \gamma_{x,F,q}$ and $\varepsilon \gamma_{x,G,q}$, respectively. These balls are defined concerning the $p$-Wasserstein distance in $\mathbb{R}$. The specific $p$-Wasserstein distance employed depends on both $\Phi$ and the support of $\xi$. We will discuss this in greater detail in the subsequent subsections. Note that the $\Phi$ of (\ref{StochsticProgWithExpectConst}) is designed for probability distributions supported on $\Xi\subseteq\mathbb{R}^{m}$, while the $\Phi$ defined in (\ref{MarkovizRobust}) is designed for probability distributions supported on $\mathbb{R}$. However, every $\Phi$ defined in (\ref{StochsticProgWithExpectConst}) induces an $\Phi$ in the context of (\ref{MarkovizRobust}). To illustrate this point, let us consider a specific case where $\Phi\left(F(x,\xi),\xi\sim\mathbb{P}\right):=\mathrm{Var}_{\xi\sim\mathbb{P}}\left[F(x,\xi)\right]$. In this situation, we have $\Phi\left(\zeta,\xi\sim\mathbb{Q}\right):=\mathrm{Var}_{\zeta\sim\mathbb{Q}}\left[\zeta\right]$. Consequently, the functions $\Phi$ examined in this work should possess this characteristic, leading our analysis to focus on cases where $\Phi$ represents a variance or an expected value.

Lastly, we adopt the letter ``A,cst'' to signify ``alternative with constraints'' in order to distinguish our proposal from the formulation presented in (\ref{MarkovizRobustGeneralWass}). From this point forward, the optimal solutions of (\ref{MarkovizRobust}) will be denoted by $\widehat{x}_{N,p,q}^{\mathrm{A,cst}}(\varepsilon)$.

\begin{figure}[t] 
 \centering
 \begin{tabular}{c}
    \includegraphics[scale=0.31]{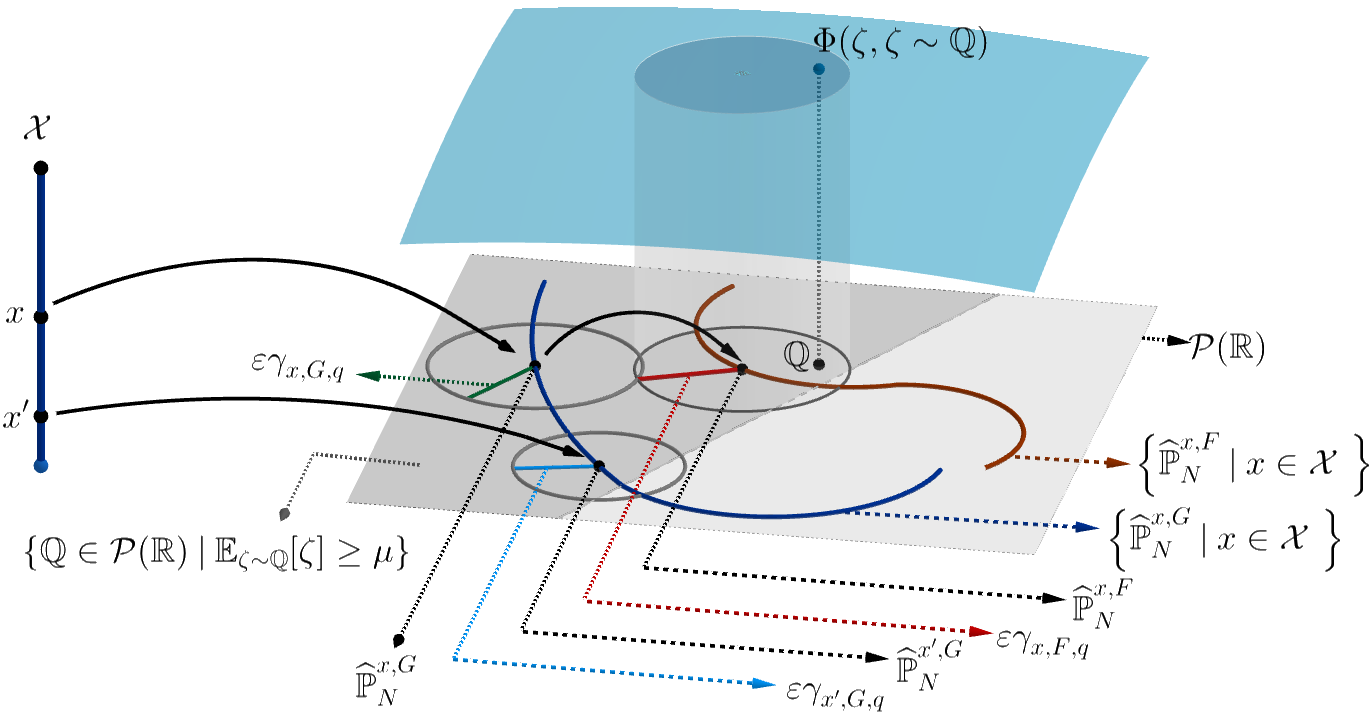}  \\
 \end{tabular}
\caption{Illustration of the optimization process proposed in (\ref{MarkovizRobust}). In this case, $x$ represents a feasible point while $x'$ represents a non-feasible one. } \label{fig:Diagrama}
\end{figure}

 Figure \ref{fig:Diagrama} shows the optimization process represented by solving  (\ref{MarkovizRobust}). This figure allows us to identify the differences between the approach proposed in this work and the standard approach represented in equation (\ref{MarkovizRobustGeneralWass}) and illustrated in Figure \ref{fig:DiagramaStandard}. Indeed, although both approaches share the fact that a given decision $x\in\mathcal{X}$ induces an optimization problem, which in both cases is a maximization problem, the difference lies in how decision $x$ influences these maximization problems. 

In Figure \ref{fig:Diagrama}, for the case of approach (\ref{MarkovizRobust}), it is observed that $x$ has influence only on the space $\mathcal{P}(\mathbb{R})$ of the probabilities supported in $\mathbb{R}$ that acts as the feasible set for the maximization problem associated with this approach. This influence consists of collapsing $\mathcal{X}$ to the space $\mathcal{P}(\mathbb{R})$ twice by assigning $x\mapsto \widehat{\mathbb{P}}_{N}^{x,G}$ and $x\mapsto \widehat{\mathbb{P}}_{N}^{x,F}$, which induces two regions in $\mathcal{P}(\mathbb{R})$. The first region (dark blue curve) is associated with the first assignment, and it determines the feasibility of the problem, as the feasible points $x$ are those that induce balls centered at $\widehat{\mathbb{P}}_{N}^{x,G}$ contained in the region $\left\{\mathbb{Q}\in\mathcal{P}(\mathbb{R})\:|\: \mathbb{E}_{\zeta\sim\mathbb{Q}}[\zeta]\geq \mu\right\}$ (dark gray region). The second region (dark red curve) is associated with the second assignment; the maximization is performed in this region only for the $x$ that are feasible, so both regions are related. 


The following subsections present the reformulation of (\ref{MarkovizRobust}) for two types of $\Phi$ functions.

\subsection{Risk neutral case}

In this scenario, we assume that $\Phi\left(F(x,\xi),\xi\sim\mathbb{P}\right):=\mathbb{E}_{\xi\sim\mathbb{P}}\left[F(x,\xi)\right]$. The initial task involves reformulating (\ref{MarkovizRobust}) as an optimization problem with finite-dimensional variables. This reformulation relies on the image of the support of $\xi$ under the functions $F(x,\cdot)$ and $G(x,\cdot)$ for each $x \in \mathcal{X}$. This goal is accomplished in the following theorem.

\begin{theorem} \label{Prop:Reformul1MarkovizDROWReformLargaLargaSupportAcot}
\hfill
\begin{enumerate}
    \item[(a)]  If $p=1$, $F(x,\Xi)$ and $G(x,\Xi)$ are intervals for each $x\in\mathcal{X}$, then  the optimization problem (\ref{MarkovizRobust})  is equivalent to the following optimization problem
\begin{equation}\label{MarkovizDROWReformLargaSupportAcot}
     \left\{\begin{array}{ll}\underset{ x\in\mathcal{X} }{\mathrm{minimize}} &   \min\left\{\frac{1}{N}{\displaystyle\sum_{i=1}^{N}}F\left(x,\widehat{\xi}_{i}\right)+\varepsilon\gamma_{x,F,q},\sup_{\xi\in\Xi}F(x,\xi)\right\} \\[0.4cm] 
     \mbox{subject to} & \max\left\{{\displaystyle \frac{1}{N}\sum_{i=1}^{N}G\left( x,\widehat{\xi}_{i}\right)-\varepsilon\gamma_{x,G,q}},\inf_{\xi\in\Xi}G(x,\xi) \right\}  \geq \mu. \end{array}\right.
\end{equation}

\item[(b)] If $p\geq 1$, $\sup_{\xi\in\Xi}F(x,\xi)=\infty$, $\inf_{\xi\in\Xi}G(x,\xi)=-\infty$, and $F(x,\Xi)$ and $G(x,\Xi)$ are intervals, all this for each  $x\in\mathcal{X}$, then the optimization problem (\ref{MarkovizRobust})  is equivalent to the following optimization problem
\begin{equation}\label{MarkovizDROWReformLargaSupportNoAcot}
     \left\{\begin{array}{ll}\underset{ x\in\mathcal{X}  }{\mathrm{minimize}} &   \frac{1}{N}{\displaystyle\sum_{i=1}^{N}}F\left(x,\widehat{\xi}_{i}\right)+\varepsilon\gamma_{x,F,q}\left(\frac{1}{p}+\frac{p-1}{p^{1/(p-1)}}\right) \\[0.4cm] 
     \mbox{subject to} & {\displaystyle \frac{1}{N}\sum_{i=1}^{N}G\left( x,\widehat{\xi}_{i}\right)-\varepsilon\gamma_{x,G,q} \left(\frac{1}{p}+\frac{p-1}{p^{1/(p-1)}}\right)}  \geq \mu. \end{array}\right.
\end{equation}
\end{enumerate}
\end{theorem}

This theorem also presents the opportunity to combine cases (i) and (ii). For instance, $F$ may satisfy the hypotheses of case (i) while $G$ satisfies the hypotheses of case (ii). In such a situation, the resulting reformulation becomes an optimization problem that incorporates the objective function of (\ref{MarkovizDROWReformLargaSupportAcot}) and the constraint of  (\ref{MarkovizDROWReformLargaSupportNoAcot}).


It is evident that the two resulting problems from Theorem \ref{Prop:Reformul1MarkovizDROWReformLargaLargaSupportAcot} are not semi-infinite, which highlights a distinct difference with the standard approach, where reformulations tend to be semi-infinite, as presented in equation (\ref{Eqn:ReformulacionDROWconConstraintRn}).

Another aspect of the approach proposed in this paper concerns feasibility. In this situation, feasibility depends on two parameters involved in the problem, specifically $\mu$ and $\varepsilon$. The subsequent corollary identifies the values of these two parameters for which problems (\ref{MarkovizDROWReformLargaSupportAcot})  and (\ref{MarkovizDROWReformLargaSupportNoAcot}) are feasible.

\begin{corollary} \label{Corol:MuAndRadioFactibleMeanVarSupportAcot}
Let $\widehat{\mu}_{N}^{\mathrm{max}}:={\displaystyle\sup_{x\in\mathcal{X}}}\frac{1}{N}\sum\limits_{i=1}^{N}G\left( x,\widehat{\xi}_{i}\right)$ and, for a given $\mu$, $ \widehat{\varepsilon}_{N,p,q}^{\mathrm{max}}(\mu):={\displaystyle\sup_{x\in\mathcal{X}}}\frac{\frac{1}{N}\sum\limits_{i=1}^{N}G(x,\widehat{\xi}_{i})-\mu}{\gamma_{x,G,q}\left(\frac{1}{p}+\frac{p-1}{p^{1/(p-1)}}\right)}$. Then the feasibility of (\ref{MarkovizDROWReformLargaSupportAcot})  and (\ref{MarkovizDROWReformLargaSupportNoAcot})  is only possible in the following cases:
\begin{enumerate}
    \item[(i)] If there exists $x\in\mathcal{X}$ such that ${\displaystyle\inf_{\xi\in\Xi} }G(x,\xi)\geq \mu$, then the optimization problem (\ref{MarkovizDROWReformLargaSupportAcot}) is feasible  for all $\varepsilon>0$.
    \item[(ii)] If ${\displaystyle\inf_{\xi\in\Xi} }G(x,\xi)< \mu$ for each $x\in\mathcal{X}$, then (\ref{MarkovizDROWReformLargaSupportAcot}) is feasible if $\mu$ and $\varepsilon$ satisfies the inequalities  $\mu < \widehat{\mu}_{N}^{\mathrm{max}}$ and $\varepsilon< \widehat{\varepsilon}_{N,1,q}^{\mathrm{max}}(\mu)$.    
    \item[(iii)]  The optimization problem (\ref{MarkovizDROWReformLargaSupportNoAcot}) is feasible if $\mu$ and $\varepsilon$ satisfies the inequalities $\mu< \widehat{\mu}_{N}^{\mathrm{max}}$ and $\varepsilon< \widehat{\varepsilon}_{N,p,q}^{\mathrm{max}}(\mu)$.  
\end{enumerate}
\end{corollary}

In this corollary we defined the expressions $\widehat{\mu}_{N}^{\mathrm{max}}$ and $\widehat{\varepsilon}_{N,p,q}^{\mathrm{max}}(\mu)$. The expression $\widehat{\mu}_{N}^{\mathrm{max}}$ represents the maximum $\mu$ level that can be targeted for the problem to remain feasible. This is consistent with what would occur in problem (\ref{StochsticProgWithExpectConst}) if the distribution $\mathbb{P}$ were known. In that problem, the phenomenon of infeasibility can also occur. However, since $\mathbb{P}$ is not known, it is impossible to determine the exact values of $\mu$ that avoid infeasibility. In this regard, $\widehat{\mu}_{N}^{\mathrm{max}}$ can be viewed as a means to identify those values. Regarding $\widehat{\varepsilon}_{N,p,q}^{\mathrm{max}}(\mu)$, provided that $\mu$ satisfies $\mu<\widehat{\mu}_{N}^{\mathrm{max}}$, this expression represents the maximum $\varepsilon$ that can be considered, ensuring the problem remains feasible. Moreover, although the corollary does not address the case where $\varepsilon= \widehat{\varepsilon}_{N,p,q}^{\mathrm{max}}(\mu)$, this case can be considered as long as an $x\in\mathcal{X}$ exists that reaches the supremum defining $\widehat{\varepsilon}_{N,p,q}^{\mathrm{max}}(\mu)$. Similarly, it is crucial to note that if $\mu=\widehat{\mu}_{N}^{\mathrm{max}}$, then $\widehat{\varepsilon}_{N,p,q}^{\mathrm{max}}(\mu)=0$ as long as an $x\in\mathcal{X}$ exists that reaches the supremum defining $\widehat{\mu}_{N}^{\mathrm{max}}$. Finally, the two expressions defined in this corollary are useful when implementing this approach in a specific situation, which will be discussed in subsequent sections.

\subsection{Variance case}

In this instance, we assume that $\Phi\left(F(x,\xi),\xi\sim\mathbb{P}\right):=\mathrm{Var}_{\xi\sim\mathbb{P}}\left[F(x,\xi)\right]$. Analogous to the previous case, the primary task is to reformulate (\ref{MarkovizRobust}) with the objective of making it computationally tractable. However, we restrict our consideration to the $2$-Wasserstein distance, since it initially allowed a reformulation that is arguably more manageable. This does not preclude the possibility of achieving a similar reformulation for other values of $p$, which is, in fact, a topic we are presently exploring.

\begin{theorem}\label{Prop:Reformul1MarkovizDROWReformLargaLarga}
Assuming $p=2$, if $F(x,\Xi)=[0,\infty)$ or $F(x,\Xi)=\mathbb{R}$, and $G(x,\Xi)=[0,\infty)$ or $G(x,\Xi)=\mathbb{R}$ for each $x\in\mathcal{X}$, then  the optimization problem (\ref{MarkovizRobust})  is equivalent to the following optimization problem with finite-dimensional variables
\begin{equation}\label{MarkovizDROWReformLarga}
     \left\{\begin{array}{ll}\underset{ x\in\mathbb{R}^{m}  }{\mathrm{minimize}} &  \left( \sqrt{\frac{1}{N}{\displaystyle\sum_{i=1}^{N}}F\left(x,\widehat{\xi}_{i}\right) ^{2}-\frac{1}{N^{2}}\left({\displaystyle\sum_{i=1}^{N}}F\left( x,\widehat{\xi}_{i}\right)\right)^{2} }+\varepsilon\gamma_{x,F,q} \right)^{2} \\[0.4cm] \mbox{subject to} & {\displaystyle \frac{1}{N}\sum_{i=1}^{N}G\left( x,\widehat{\xi}_{i}\right)-\varepsilon\gamma_{x,G,q} \geq \mu,} \\[0.4cm] 
&  x\in\mathcal{X} \end{array}\right.
\end{equation}
\end{theorem}

It is noteworthy that, thus far, no reformulation for the standard version (\ref{MarkovizRobustGeneralWass}) exists for this case. This is because the consideration of variance as an objective function from a distributionally robust perspective has not been widely explored in a general context. Some studies have focused on specific instances, such as \cite{Blanchet2018}, which direct their analysis toward the context of portfolio optimization. In contrast, the approach proposed in this paper permits the consideration of various types of functions $F$. Furthermore, if $F$, $\gamma_{x,F,q}$, and $\gamma_{x,F,q}$ are convex with respect to $x$, then the optimization problem (\ref{MarkovizDROWReformLarga}) is convex, which may prove advantageous during implementation.

Lastly, it is worth addressing the matter of feasibility in this case. In this regard, the constraints of (\ref{MarkovizDROWReformLarga}) coincide with those of (\ref{MarkovizDROWReformLargaSupportNoAcot}), which implies that the values of $\mu$ and $\varepsilon$ for which(\ref{MarkovizDROWReformLarga}) is feasible correspond to those established in part (iii) of Corollary \ref{Corol:MuAndRadioFactibleMeanVarSupportAcot}. The proofs for Theorems \ref{Prop:Reformul1MarkovizDROWReformLargaLargaSupportAcot} and \ref{Prop:Reformul1MarkovizDROWReformLargaLarga} , as well as Corollary \ref{Corol:MuAndRadioFactibleMeanVarSupportAcot}, can be found in Appendix \ref{Sec:Appendix:ProofThmMeanVarianze}.

\subsection{Expected confidence level of Wasserstein radius} \label{Subsec:Bootstrap}

To conclude this section, we aim to estimate the probability that solutions generated by the proposed approach satisfy the constraint of (\ref{StochsticProgWithExpectConst}) when evaluated out-of-sample. Obtaining an estimate of this probability may suggest a method for choosing $\varepsilon$. For instance, in (\ref{MarkovizRobust}), if the priority were to select $\varepsilon$ such that $\mathbb{E}_{\xi\sim\mathbb{P}}\left[G\left( \widehat{x}_{N,p,q}^{A,cst}(\varepsilon),\xi\right)\right] \geq \mu$ is satisfied with high probability, where $\widehat{x}_{N,p,q}^{A,cst}(\varepsilon)$ is an optimal solution of (\ref{MarkovizRobust}), then the following lemma could offer a criterion to accomplish this.

\begin{lemma} \label{Lemma:Confianza}
Let $\varepsilon>0$ such that $\mathbb{P}\in\mathcal{B}_{\varepsilon}(\widehat{\mathbb{P}}_{N})$. If $\varepsilon$ is such that (\ref{MarkovizRobust}) is feasible, and  $\widehat{x}_{N,p,q}^{A,cst}(\varepsilon)$ is a optimal solution of (\ref{MarkovizRobust}), then  $\mathbb{E}_{\xi\sim\mathbb{P}}\left[G\left( \widehat{x}_{N,p,q}^{A,cst}(\varepsilon),\xi \right)\right] \geq \mu$.
\end{lemma}
The lemma suggests that finding an $\varepsilon$ such that $\mathbb{P}\in\mathcal{B}_{\varepsilon}(\widehat{\mathbb{P}}_{N})$ is sufficient, and this can be achieved with large values of $\varepsilon$. However, it is essential to consider that $\varepsilon$ might be constrained by $\widehat{\varepsilon}_{N,p,q}^{\mathrm{max}}(\mu)$ in (\ref{MarkovizRobust}). Another challenge lies in the absence of an efficient method for determining the value of $\varepsilon$ from which the condition $\mathbb{P}\in\mathcal{B}_{\varepsilon}(\widehat{\mathbb{P}}_{N})$ is satisfied.

With this in mind, our approach does not focus on ensuring the satisfaction of this condition. Instead, we concentrate on estimating the probability that $\mathbb{E}_{\xi\sim\mathbb{P}}\left[G\left( \widehat{x}_{N}(\varepsilon),\xi \right)\right] \geq \mu$ is satisfied, with this probability taken with respect to the sample. To accomplish this, we employ a strategy based on the Bootstrap method, assuming that $\xi$ meets all the necessary conditions for the Bootstrap-based technique to be valid (see \cite{efron1994}). Given $\varepsilon>0$, our objective is to estimate the probability that the constraint is satisfied, which we term the \textit{expected confidence level of $\varepsilon$}. The following method provides an estimate of this probability.
\begin{enumerate}
\item[$\bullet$] Given $\varepsilon>0$, we generate $K$ bootstrap samples of size $N$ with repetition from $\widehat{\Xi}_{N}=\left\{\widehat{\xi}_{1},\ldots,\widehat{\xi}_{N}\right\}$. These bootstrap samples are denoted as $\widehat{\Xi}_{N,i}^{\mathrm{bt}}:=\left\{\widehat{\xi}_{1}^{\mathrm{bt},i},\ldots,\widehat{\xi}_{N}^{\mathrm{bt},i}\right\}$ for $i=1,2\ldots,K$.\\

We then partition each $\widehat{\Xi}_{N,i}^{\mathrm{bt}}$ into a training dataset $\widehat{\Xi}_{N_{T},i}^{\mathrm{bt}}$ of size $N_{T}$ and a validation dataset $\widehat{\Xi}_{N_{V},i}^{\mathrm{bt}}$ of size $N_{V}=N-N_{T}$. For each training dataset $\widehat{\Xi}_{N_{T},i}^{\mathrm{bt}}$, we compute $\widehat{x}_{N_{T},p,q}^{A,cst}(\varepsilon)$, denoting it as $\widehat{x}_{N_{T},p,q}^{A,cst(i)}(\varepsilon)$ to emphasize its association with the sample $\widehat{\Xi}_{N_{T},i}^{\mathrm{bt}}$.\\

Additionally, we define $\widehat{\mathbb{P}}_{N_{V},i}^{\mathrm{bt}}$ as the empirical distribution generated by the validation dataset $\widehat{\Xi}_{N_{V},i}^{\mathrm{bt}}$. Subsequently, we calculate the sample mean of $G\left( \widehat{x}_{N_{T},p,q}^{A,cst(i)}(\varepsilon),\xi \right)$ induced by the validation dataset $\widehat{\mathbb{P}}_{N_{V},i}^{\mathrm{bt}}$, denoted as $\mathbb{E}_{\xi\sim\widehat{\mathbb{P}}_{N_{V},i}^{\mathrm{bt}}}\left[G\left(\widehat{x}_{N_{T},p,q}^{A,cst(i)}(\varepsilon),\xi \right)\right]$. Consequently, the expected confidence level of $\varepsilon$ is the percentage of occurrences in which $\mathbb{E}_{\xi\sim\widehat{\mathbb{P}}_{N_{V},i}^{\mathrm{bt}}}\left[G\left( \widehat{x}_{N_{T},p,q}^{A,cst(i)}(\varepsilon),\xi \right)\right]\geq \mu$.
\end{enumerate}
If the estimated confidence level of $\varepsilon$ is $\beta$, it implies that the constraint is satisfied with an approximate probability of $\beta/100$. Therefore, if the priority is to satisfy the constraint, this concept of expected confidence level could suggest a method for choosing $\varepsilon$. The strategy might involve finding the smallest $\varepsilon$ for which the expected confidence level $\beta$ is acceptable. As the confidence level decreases when $\varepsilon$ becomes smaller, the maximum expected confidence level attainable corresponds to $\varepsilon= \widehat{\varepsilon}_{N,p,q}^{\mathrm{max}}(\mu)$. Thus, the search for an $\varepsilon$ meeting the required expected confidence level can commence from this value when it is finite. In Section \ref{Sec:Experimentsesults}, we will explore situations where the expected confidence level is calculated and examine the range of values it can assume. The proof of Lemma \ref{Lemma:Confianza} is relegated to \ref{Apendice:PruebasLema}.

Finally, it is possible that, for every feasible $\varepsilon$, the expected confidence level is 100\%, meaning that all decisions induced by the approach proposed for each $\varepsilon$ satisfy the constraint with high probability. In such cases, if the priority was to satisfy the constraint, the focus should shift to selecting an $\varepsilon$ that induces the lowest out-of-sample value of the objective function $\mathbb{E}_{\xi\sim\mathbb{P}}\left[F\left(\widehat{x}_{N,p,q}^{A,cst}(\varepsilon),\xi\right)\right]$. Since the true probability distribution $\mathbb{P}$ is unknown, it is recommended to use bootstrap-based methods such as those presented in \cite{Kuhn2018} to estimate this expression. These methods are commonly known as the Holdout method and the $k$-fold cross-validation method.

\section{Portfolio optimization}  \label{Sec:PortfolioSubsection}

In this section, we direct our analysis toward the context of portfolio optimization, setting the stage for the presentation of numerical results. Specifically, we focus on two problems. The first involves the objective function in (\ref{StochsticProgWithExpectConst}) being the Conditional Value at Risk (CVaR). The second concerns the objective function in (\ref{StochsticProgWithExpectConst}) being the variance, known as the mean-variance problem. Both problems fall within the two types of $\Phi$ that we aimed to analyze in this paper.

Before delving into the details, it is essential to describe the portfolio optimization context. In this setting, we have $m=n$, where $m$ represents the number of assets, $\xi_i$ denotes the return of the $i$-th asset, and $x_{i}$ signifies the proportion of the initial amount invested in the $i$-th asset. As the returns of each asset are random with an \emph{unknown} distribution $\mathbb{P}$, $\xi=(\xi_{1},\ldots,\xi_{m})\in\mathbb{R}^{m}$ is a random vector. Furthermore, $x=(x_{1},\ldots,x_{m})\in\mathbb{R}^{m}$ is a portfolio, which is a vector of weights satisfying the relation $\sum_{i=1}^{m}x_{i}=1$ and additional convex constraints acceptable to the investor, described by the set $\mathsf{X}$. Consequently, we have $\mathcal{X}=\left\{x\in\mathbb{R}^{m}\: :\:\sum_{i=1}^{m} x_{i}= 1,\: x\in\mathsf{X}\right\}$. Moreover, $G(x,\xi)=\langle x,\xi\rangle$, where $\langle\cdot,\cdot\rangle$ denotes the Euclidean inner product in $\mathbb{R}^{m}$. Then, $\langle x,\xi\rangle$ represents the return of the portfolio $x$. Additionally, we will focus our analysis on the case $q=p=2$. In this regard, note that $G$ satisfies Assumption \ref{AssumptionPrincipal} because $G$ is a $2$-Lipschitz functions with respect to $\xi$ with Lipschitz constant $\gamma_{x,G,q}=\|x\|_{2}$. $F$ is specified in each of the cases we will analyze. To simplify notation, we will write $\|\cdot\|$ instead of $\|\cdot\|_{2}$. Lastly, $\mu$ is the minimum level of return acceptable to the investor.

\subsection{Mean-Risk Portfolio Optimization}

In this part, we concentrate on the case where the objective function is the Conditional Value at Risk (CVaR). In such problems, the objective is to find a portfolio that minimizes the CVaR, subject to the expected return induced by the portfolio exceeding a specified level. The formulation of this problem is as follows:
\begin{equation}\label{StochsticProgWithExpectConsCVaR}
J^{\mathrm{CVaR}_{\mathrm{cst}}}:=\left\{\begin{array}{ll}
{\displaystyle \min_{x\in\mathbb{R}^{m}} } & \mathrm{CVaR}_{\alpha,\xi\sim\mathbb{P}}\left(-\langle x,\xi\rangle\right)\\
\mbox{subject to} & \mathbb{E}_{\xi\sim\mathbb{P}}\left[\langle x,\xi\rangle\right]\geq \mu,\\
& x\in\mathcal{X},
\end{array}  \right.
\end{equation}
Here, the subscript ``cst'' indicates that the problem is subject to a stochastic constraint. Considering that the CVaR can be expressed as the minimization of an expected value \cite{Rockafellar2000}, it follows that (\ref{StochsticProgWithExpectConsCVaR}) can be reformulated as:
\begin{equation}\label{StochsticProgWithExpectConsCVaRReform}
\left\{\begin{array}{ll}
{\displaystyle \min_{x\in\mathbb{R}^{m},\tau\in\mathbb{R}} } & \mathbb {E}_{\xi\sim\mathbb{P}}\left[\max\left\{-\frac{\langle x,\xi\rangle}{\alpha}+\left(1-\frac{1}{\alpha}\right)\tau,\tau\right\}\right]\\
\mbox{subject to} & \mathbb{E}_{\xi\sim\mathbb{P}}\left[\langle x,\xi\rangle\right]\geq \mu,\\
& x\in\mathcal{X},
\end{array}  \right.
\end{equation}
\sloppy This problem seeks to minimize the expected value of losses within the $\alpha$\% highest losses, subject to the constraint on the expected value of returns. In this context, we consider $F(x,\tau,\xi)=\max\left\{-\frac{\langle x,\xi\rangle}{\alpha}+\left(1-\frac{1}{\alpha}\right)\tau,\tau\right\}$, where the decision variables are $x$ and $\tau$. This function $F$ is $2$-Lipschitz with a Lipschitz constant  $\gamma_{x,F,q}=\frac{\|x\|}{\alpha}$. Consequently, Theorem \ref{Prop:Reformul1MarkovizDROWReformLargaLargaSupportAcot} in part (b) enables the reformulation of the distributionally robust version (\ref{MarkovizRobust}) for this case as:
\begin{equation}\label{MarkovizDROWReformLargaCVaR}
     \widehat{J}_{N}^{\mathrm{CVaR}_{\mathrm{cst}}}(\varepsilon):=\left\{\begin{array}{ll}\underset{ x\in\mathbb{R}^{m}, \tau\in\mathbb{R}  }{\mathrm{minimize}} &  \frac{\varepsilon}{\alpha}\|x\| +\frac{1}{N}\sum\limits_{i=1}^{N}\max\left\{-\frac{\langle x,\widehat{\xi}_{i}\rangle}{\alpha}+\left(1-\frac{1}{\alpha}\right)\tau,\tau\right\} \\[0.4cm] \mbox{subject to} & {\displaystyle \frac{1}{N}\sum_{i=1}^{N}\langle x,\widehat{\xi}_{i}\rangle-\varepsilon\|x\| \geq \mu,} \\[0.4cm] 
&  x\in\mathcal{X}. \end{array}\right.
\end{equation}
By introducing auxiliary variables, this problem can be expressed as a second-order cone program, resulting in a convex problem. In contrast, applying the standard approach presented in (\ref{MarkovizRobustGeneralWass}) to this case yields, through (\ref{Eqn:ReformulacionDROWconConstraintRn}) and some algebraic development, the following problem:
\begin{equation}\label{MarkovizDROWReformLargaCVaRStandar}
     \left\{\begin{array}{ll}\underset{ x\in\mathbb{R}^{m}, \tau\in\mathbb{R}, \lambda> 0  }{\mathrm{minimize}} &  \lambda \varepsilon^{2}+\frac{1}{N}\sum\limits_{i=1}^{N}\max\left\{ \frac{\|x\|^{2}}{4\lambda \alpha^{2}} -\frac{\langle x,\widehat{\xi}_{i}\rangle}{\alpha}+\left(1-\frac{1}{\alpha}\right)\tau,\tau\right\} \\[0.4cm] \mbox{subject to} & {\displaystyle \frac{1}{N}\sum_{i=1}^{N}\langle x,\widehat{\xi}_{i}\rangle-\varepsilon\|x\| \geq \mu,} \\[0.4cm] 
&  x\in\mathcal{X}. \end{array}\right.
\end{equation}
From this, it is evident that the standard approach produces a problem with more variables, and one of those variables appears as a divisor, complicating the optimization process. Furthermore, with some algebraic development, this problem could be recast as a program with quadratic and second-order cone constraints. However, the matrices describing the quadratic constraints are undefined, potentially rendering this problem non-convex. Therefore, this case suggests that the approach proposed in this paper is a viable alternative to consider when addressing a stochastic problem with expected value constraints from a distributionally robust perspective.

\subsection{Mean-variance Portfolio Optimization}

In this part, we examine the Markowitz mean-variance portfolio selection optimization problem. The objective is to choose portfolio weights that minimize the variance of the return rate while adhering to a constraint on the expected value of the return rate. The formulation of this problem is as follows:
\begin{equation}\label{StochsticProgWithExpectConstVariance}
J^{\mathrm{Var}_{\mathrm{cst}}}:=\left\{\begin{array}{ll}
{\displaystyle \min_{x\in\mathbb{R}^{m}} } & \mathbb{V}\mathrm{ar}_{\xi\sim\mathbb{P}}\left[\langle x,\xi\rangle\right]\\
\mbox{subject to} & \mathbb{E}_{\xi\sim\mathbb{P}}\left[\langle x,\xi\rangle\right]\geq \mu,\\
& x\in\mathcal{X},
\end{array}  \right.
\end{equation}
The notation ``$\mathrm{Var}_{\mathrm{cst}}$'' indicates that the problem aims to minimize the variance and incorporates a stochastic constraint. Additionally, in this case, we are considering $F=G$. Bearing this in mind, Theorem \ref{Prop:Reformul1MarkovizDROWReformLargaLarga} allows us to deduce that the distributionally robust version (\ref{MarkovizRobust}) for this case can be reformulated as:
\begin{equation}\label{MarkovizDROWReformLargaVariance}
     \left\{\begin{array}{ll}\underset{ x\in\mathbb{R}^{m}  }{\mathrm{minimize}} &  \left( \sqrt{\frac{1}{N}{\displaystyle\sum_{i=1}^{N}}\langle x,\widehat{\xi}_{i}\rangle ^{2}-\frac{1}{N^{2}}\left({\displaystyle\sum_{i=1}^{N}}\langle x,\widehat{\xi}_{i}\rangle\right)^{2} }+\varepsilon\|x\| \right)^{2} \\[0.4cm] \mbox{subject to} & {\displaystyle \frac{1}{N}\sum_{i=1}^{N}\langle x,\widehat{\xi}_{i}\rangle-\varepsilon\|x\| \geq \mu,} \\[0.4cm] 
&  x\in\mathcal{X} \end{array}\right.
\end{equation}
To facilitate comprehension, we can rewrite this problem in a more intuitive manner. Let $\widehat{\Sigma}_{N}$ represent the biased covariance matrix and $\widehat{\mathbf{m}}_{N}$ denote the vector of sample means, both derived from the sample $\widehat{\xi}_{1},\ldots,\widehat{\xi}_{N}$. Consequently, the aforementioned problem can be expressed as:
\begin{equation} \label{MarkovizDROWReformDeFormLargaPortflio}
\widehat{J}_{N}^{\mathrm{Var}_{\mathrm{cst}}}(\varepsilon):=\left\{\begin{array}{ll}
{\displaystyle \inf_{x\in\mathbb{R}^{m} } } & \left( \sqrt{\langle x,\widehat{\Sigma}_{N}, x\rangle} +\varepsilon\left\|x\right\| \right)^{2}\\
\mbox{subject to} & \langle \widehat{\mathbf{m}}_{N}, x\rangle-\varepsilon\left\|x\right\|\geq \mu, \\
& x\in\mathcal{X}. 
\end{array}
\right.
\end{equation}
By neglecting the square term in the objective function—which does not impact the optimal solutions—the problem can be classified as a second-order cone program. This convex optimization problem is generally not computationally expensive. Additionally, the resulting problem bears resemblance to the optimization problem that would arise from applying the Sample Average Approximation (SAA) approach. The distinction lies in the fact that (\ref{MarkovizDROWReformDeFormLargaPortflio}) includes $\varepsilon\left\|x\right\|$ in both the objective function and the constraint. This observation implies that $\varepsilon\left\|x\right\|$ serves as a form of regularization, modifying the standard SAA approach.

Additionally, it is noteworthy that reference \cite{Blanchet2018} examines problem (\ref{StochsticProgWithExpectConstVariance}) from the perspective of the distributionally robust standard approach (\ref{MarkovizRobustGeneralWass}). The resulting problem in \cite{Blanchet2018} is the same as (\ref{MarkovizDROWReformDeFormLargaPortflio}), although it originates from a different viewpoint. This observation suggests that there may be instances in which the standard approach and the proposed approach are equivalent. Nevertheless, this equivalence does not generally hold. For a more comprehensive discussion on this topic, we recommend consulting \cite{fonsecaDecDep2023}, where conditions that allow for equivalence in the case of stochastic problems without stochastic constraints are established. These conditions are also applicable to cases that involve stochastic constraints.

\subsection{Numerical experiments and results} \label{Sec:Experimentsesults}

This subsection concentrates on evaluating the performance of the proposed approach for cases where the objective function is either variance or CVaR. The analysis is conducted jointly for both objective functions. Consequently, this subsection is divided into two parts. In the first part, synthetically generated data is employed, which refers to data generated by a known distribution. In the second part, real data from the financial market is used. For both cases, we consider $\mathsf{X}=\mathbb{R}^{m}_{+}$, indicating that short selling is not taken into account.

\subsubsection{Using simulated data.}

We consider a market of $m = 10$ assets with returns following the form adopted in \cite{Kuhn2018}, where we assume $\xi_{i}=\psi +\zeta_{i}$. Here, $\psi$ and $\zeta_{i}$ are independent, with $\psi\sim \mathcal{N}(0,2\%)$ and $\zeta_{i}\sim \mathcal{N}(i\times3\% ,i\times 2.5\%)$ for each $i=1,2,\ldots,m$. Under this assumption, the assets are ordered from the lowest to the highest return and volatility. Furthermore, we denote the vector of means as $\mathbf{m}$ and the covariance matrix of $\xi$ as $\Sigma$. In this case, both $\mathbf{m}$ and $\Sigma$ can be readily calculated from the distribution of $\xi$. Given $x\in\mathbb{R}^{m}$, we define $R(x):=\mathbf{m}^{T}x$ as the expected return induced by $x$, and $V(x):=x^{T}\Sigma x$ as the variance induced by $x$ to simplify the notation. Since we possess complete information about the returns, the optimal portfolios $x^{\mathrm{Var}{\mathrm{cst}}*}$ and $x^{\mathrm{CVaR}{\mathrm{cst}}*}$ of (\ref{StochsticProgWithExpectConstVariance}) and (\ref{StochsticProgWithExpectConsCVaR}) respectively are known. Finally, as our analysis focuses on the case $p=q=2$, we will omit the influence of $p$ and $q$ in the notation of the solutions. In that sense, we denote the solutions of (\ref{MarkovizDROWReformLargaCVaR}) and (\ref{MarkovizDROWReformDeFormLargaPortflio}) as $\widehat{x}_{N}^{\mathrm{CVaR}_{\mathrm{cst}}}(\varepsilon)$ and $\widehat{x}_{N}^{\mathrm{Var}{\mathrm{cst}}}(\varepsilon)$ respectively, and $\widehat{\varepsilon}_{N}^{\mathrm{max}}(\mu)$ to refer to $\widehat{\varepsilon}_{N,p,q}^{\mathrm{max}}(\mu)$.

\paragraph{\textbf{Impact of the Wasserstein Radius $\varepsilon$.}}

\begin{figure}[tbp] 
 \centering
 \begin{tabular}{ccc}
    \includegraphics[scale=0.32]{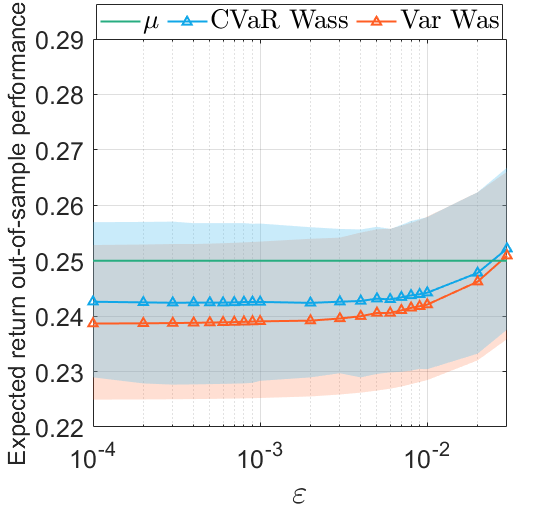} &  \includegraphics[scale=0.32]{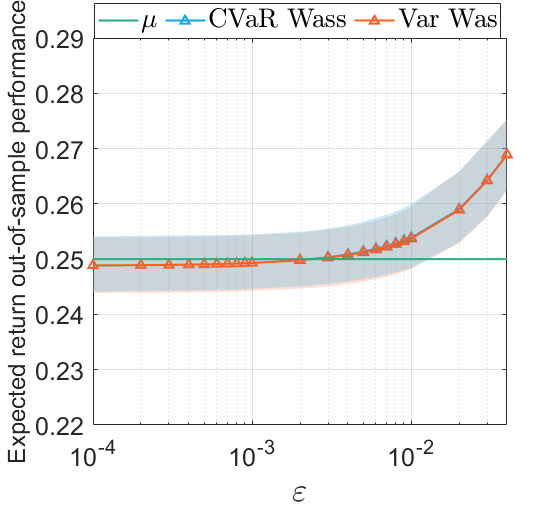} & \includegraphics[scale=0.32]{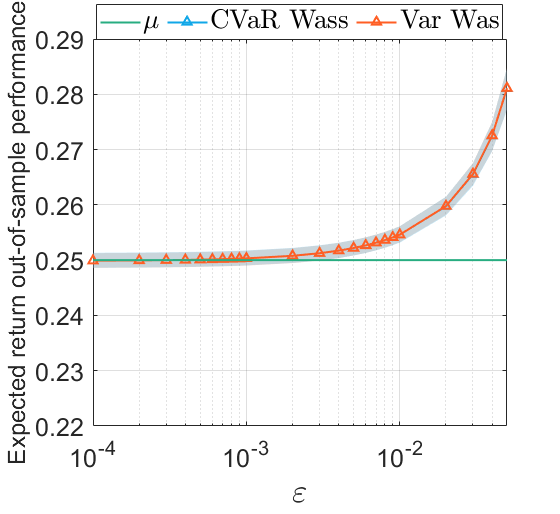} \\
      (a) $N=30.$ & (b) $N=300.$ &  (c) $N=3000.$
 \end{tabular}
\caption{The out-of-sample performance of expected returns $R\left(\widehat{x}_{N}^{\mathrm{CVaR}_{\mathrm{cst}}}(\varepsilon)\right)$ (blue) and $R\left(\widehat{x}_{N}^{\mathrm{Var}_{\mathrm{cst}}}(\varepsilon)\right)$ (orange) are illustrated as functions of the Wasserstein radius $\varepsilon$, based on 200 simulations. Solid lines represent the means, while the shaded areas correspond to the region between the 20\% and 80\% quantiles of data generated by the 200 simulations. In this case, $\mu=0.25$ and $\alpha=0.05$. } \label{fig:RetursnVsEpsilon}
\end{figure}

\begin{figure}[tbp!] 
 \centering
 \begin{tabular}{ccc}
    \includegraphics[scale=0.32]{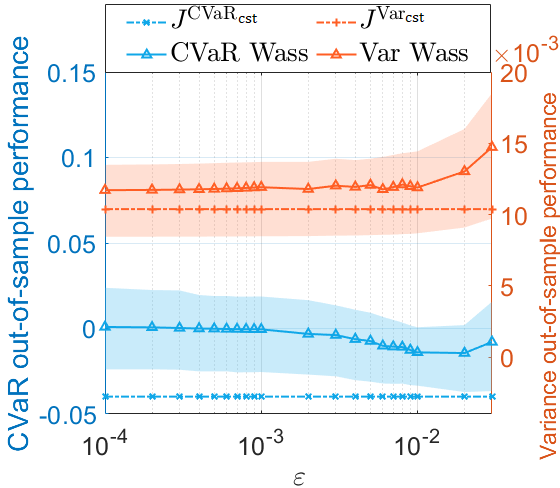} &  \includegraphics[scale=0.32]{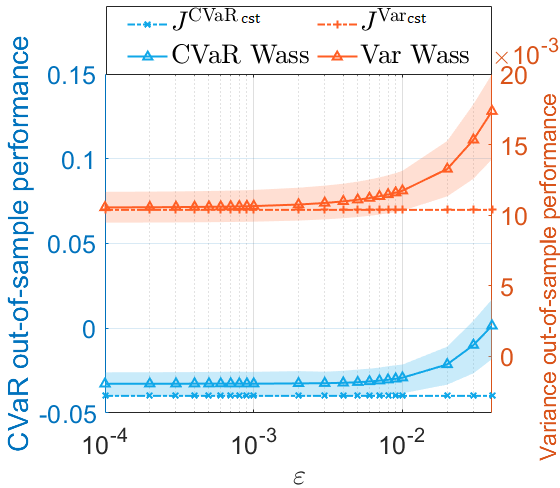} & \includegraphics[scale=0.32]{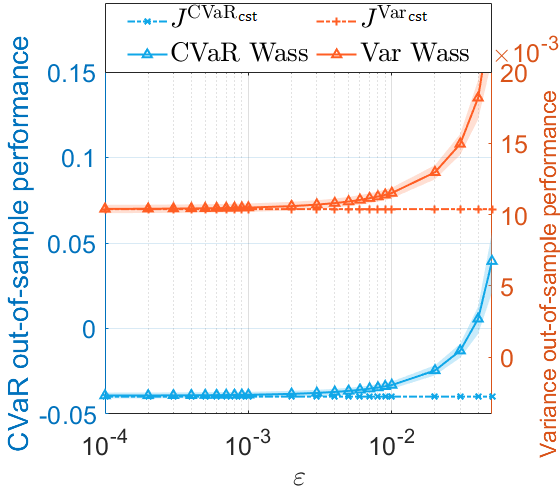} \\
     (a) $N=30.$ & (b) $N=300.$ & (c) $N=3000.$\\
     \includegraphics[scale=0.32]{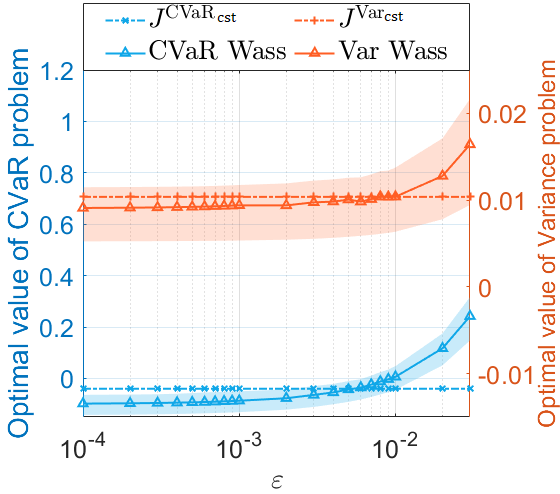} &  \includegraphics[scale=0.32]{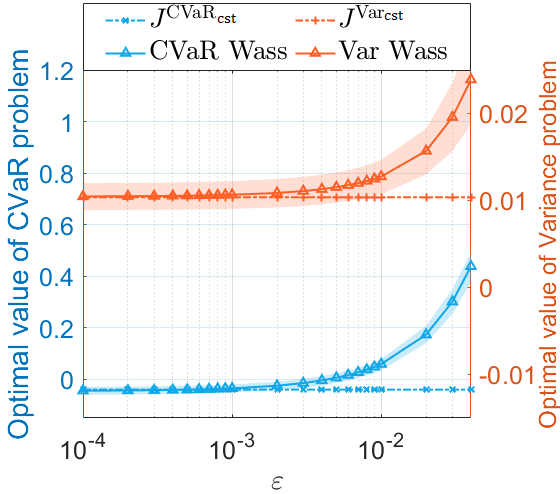} & \includegraphics[scale=0.32]{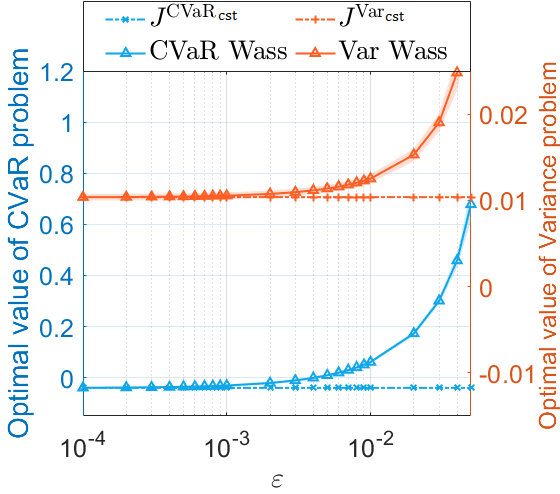} \\
     (d) $N=30.$ & (e) $N=300.$ & (f) $N=3000.$ \\
    \includegraphics[scale=0.32]{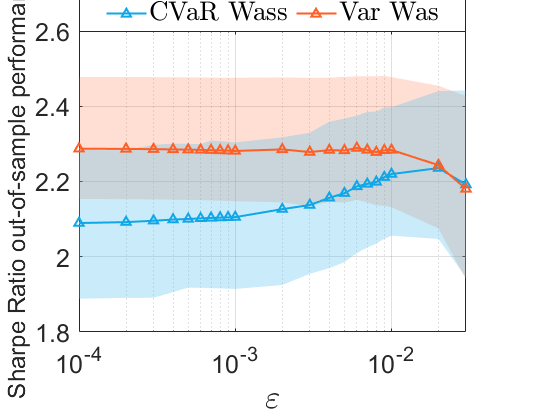} &  \includegraphics[scale=0.32]{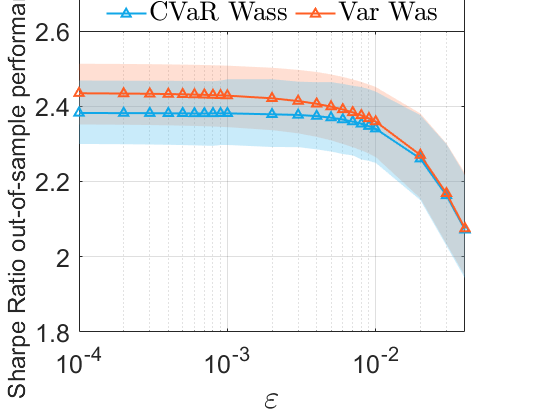} & \includegraphics[scale=0.32]{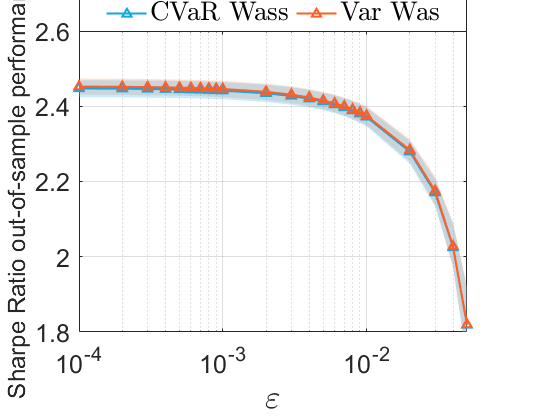} \\
      (g) $N=30.$ & (h) $N=300.$ & (i) $N=3000.$
 \end{tabular}
\caption{
This figure displays the out-of-sample performance of the variance $V\left(\widehat{x}_{N}^{\mathrm{Var}_{\mathrm{cst}}}(\varepsilon)\right)$ (orange) and the Conditional Value at Risk $\mathrm{CVaR}_{\alpha,\xi\sim\mathbb{P}}\left(-\left\langle \widehat{x}_{N}^{\mathrm{CVaR}_{\mathrm{cst}}}(\varepsilon),\xi\right\rangle\right)$ (blue), optimal values $\widehat{J}_{N}^{\mathrm{Var}_{\mathrm{cst}}}(\varepsilon)$ (orange) and $\widehat{J}_{N}^{\mathrm{CVaR}_{\mathrm{cst}}}(\varepsilon)$ (blue), and Sharpe Ratios as functions of the Wasserstein radius $\varepsilon$, estimated based on 200 simulations. The solid lines represent the means, and the shaded areas indicate the tubes between the 20\% and 80\% quantiles of data generated by 200 simulations. In this case, $\mu=0.25$ and $\alpha=0.05$.
} \label{fig:RetursnVarianceOptValueVsEpsilon}
\end{figure}

The primary aim is to investigate the effect of the Wasserstein radius $\varepsilon$ on the optimal distributionally robust portfolios generated by our approach and their out-of-sample performance. As such, we solve problems (\ref{MarkovizDROWReformLargaCVaR}) and (\ref{MarkovizDROWReformDeFormLargaPortflio}) using samples of size $N \in \left\{30, 300, 3000\right\}$. Figure \ref{fig:RetursnVsEpsilon} displays the region between the 20\% and 80\% quantiles (shaded areas) and the mean (solid lines) of the out-of-sample performance of expected returns $R\left(\widehat{x}_{N}^{\mathrm{CVaR}_{\mathrm{cst}}}(\varepsilon)\right)$ (blue) and $R\left(\widehat{x}_{N}^{\mathrm{Var}_{\mathrm{cst}}}(\varepsilon)\right)$ (orange) as functions of $\varepsilon$, estimated using 200 independent simulation runs. The same color convention is maintained for the remaining figures. It is observed that the out-of-sample performance of both expected returns increases as the Wasserstein radius $\varepsilon$ takes larger values. This observation supports the idea that larger values of $\varepsilon$ suggest a higher probability of constraint satisfaction. However, larger values of $\varepsilon$ also increase the out-of-sample values of variance $V\left(\widehat{x}_{N}^{\mathrm{Var}_{\mathrm{cst}}}(\varepsilon)\right)$ and Conditional Value at Risk $\mathrm{CVaR}_{\alpha,\xi\sim\mathbb{P}}\left(-\left\langle \widehat{x}_{N}^{\mathrm{CVaR}_{\mathrm{cst}}}(\varepsilon),\xi\right\rangle\right)$, as well as optimal values $\widehat{J}_{N}^{\mathrm{Var}_{\mathrm{cst}}}(\varepsilon)$ and $\widehat{J}_{N}^{\mathrm{CVaR}_{\mathrm{cst}}}(\varepsilon)$, as shown in Figures \ref{fig:RetursnVarianceOptValueVsEpsilon}. This encourages the preference for $\varepsilon$ values that provide a balance between the probability of satisfying the constraint and the out-of-sample value of the objective function.

This figure also reveals that Sharpe ratios $R\left(\widehat{x}_{N}^{\mathrm{CVaR}_{\mathrm{cst}}}(\varepsilon)\right)/\sqrt{V\left(\widehat{x}_{N}^{\mathrm{CVaR}_{\mathrm{cst}}}(\varepsilon)\right)}$ and $R\left(\widehat{x}_{N}^{\mathrm{Var}_{\mathrm{cst}}}(\varepsilon)\right)/\sqrt{V\left(\widehat{x}_{N}^{\mathrm{Var}_{\mathrm{cst}}}(\varepsilon)\right)}$ tend to decrease as $\varepsilon$ rises. However, it is crucial to keep in mind that large values of $\varepsilon$ should not be the preferred choices, and the most appropriate value of $\varepsilon$ is the one at which the constraint starts to be satisfied with high probability. For example, in Figure \ref{fig:RetursnVsEpsilon}(b), this situation occurs around $\varepsilon=10^{-2}$. At that value of $\varepsilon$ in Figure \ref{fig:RetursnVarianceOptValueVsEpsilon}(h), it can be observed that the Sharpe Ratio of the portfolios generated by both problems is not small. One could argue that it is acceptable since it does not represent a significant decrease compared to the Sharpe Ratio obtained by the SAA approach, which occurs when $\varepsilon=0$.

\begin{figure}[t] 
 \centering
 \begin{tabular}{ccc}
     \includegraphics[scale=0.27]{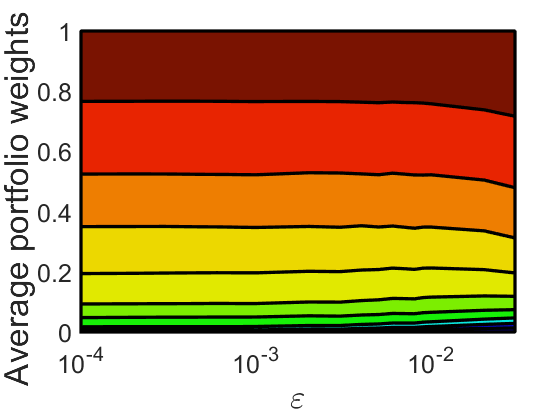} &  \includegraphics[scale=0.27]{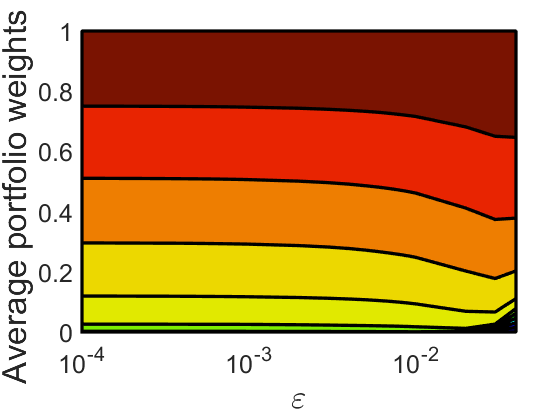} & \includegraphics[scale=0.27]{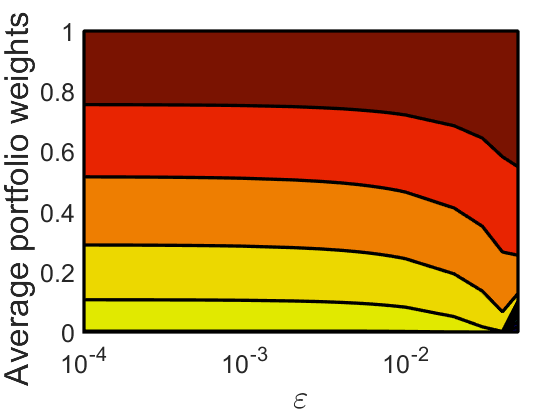} \\
     $N=30.$ & $N=300.$ & $N=3000.$ \\
     \includegraphics[scale=0.27]{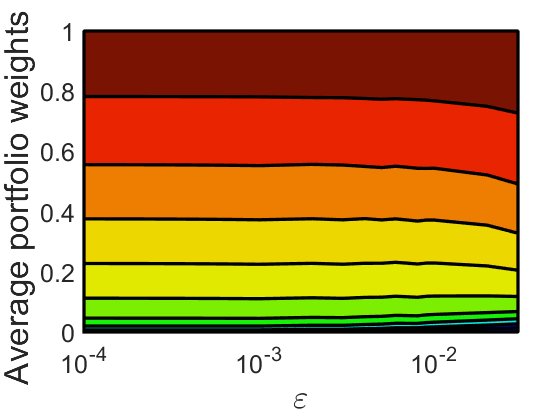} &  \includegraphics[scale=0.27]{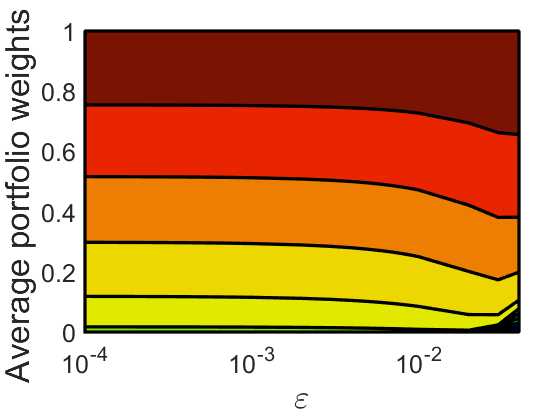} & \includegraphics[scale=0.27]{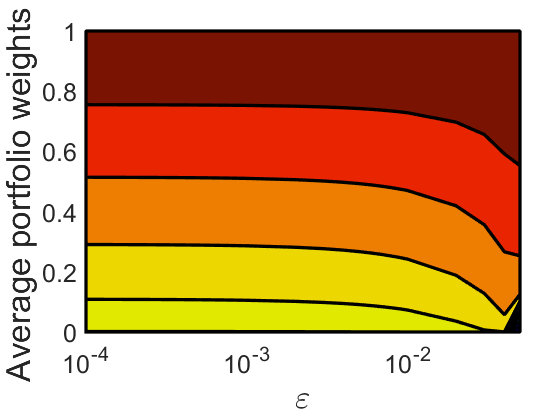} \\
     $N=30.$ & $N=300.$ & $N=3000.$
 \end{tabular}
\caption{Optimal portfolio composition as a function of the Wasserstein radius $\varepsilon$ averaged over 200 simulations. The portfolio weights are depicted in ascending order, with the weight of asset 1 at the bottom and that of asset 10 at the top. In this instance, $\mu=0.25$ and $\alpha=0.05$. } \label{fig:AverageWeightsVsEpsilon}
\end{figure}

The final aspect to examine is the behavior of the portfolios generated by the proposed approach. Figure \ref{fig:AverageWeightsVsEpsilon} illustrates the corresponding optimal portfolio weights  as a function of $\varepsilon$, averaged over 200 independent simulation runs. Our numerical results indicate that the optimal distributionally robust portfolios generated by our approach tend to assign low weight to assets with low return, even if they exhibit low volatility while allocating more weight to assets with high return, even if they possess high volatility. This occurs as the Wasserstein radius $\varepsilon$ increases. Furthermore, this figure shows that in both cases, with variance and CVaR as the objective function, the resulting portfolios exhibit similar behavior on average.

\begin{figure}[tbp] 
 \centering
 \begin{tabular}{ccc}
    \includegraphics[scale=0.32]{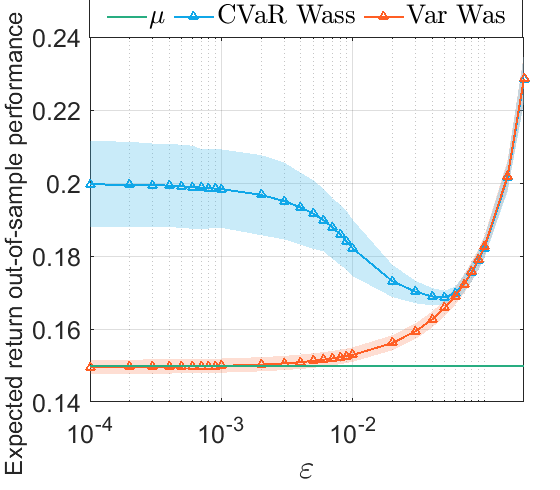} &\includegraphics[scale=0.32]{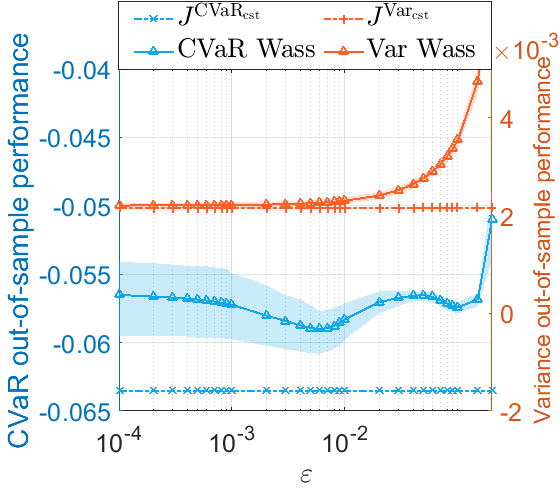} &  \includegraphics[scale=0.32]{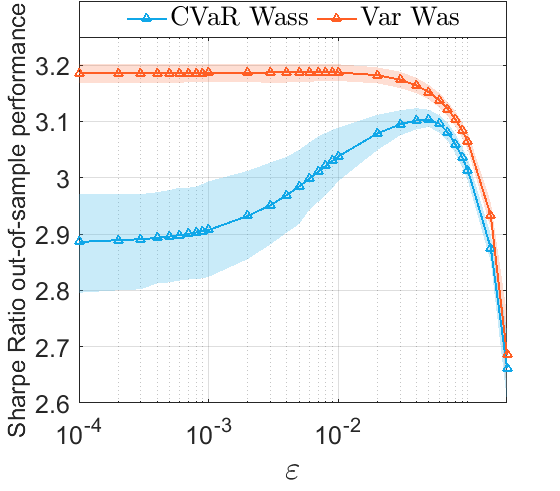}  \\
     (a)  & (b)  & (c)
 \end{tabular}
\caption{ This figure presents the out-of-sample performance of expected returns 
$R\left(\widehat{x}_{N}^{\mathrm{Var}_{\mathrm{cst}}}(\varepsilon)\right)$ (orange) and $R\left(\widehat{x}_{N}^{\mathrm{CVaR}_{\mathrm{cst}}}(\varepsilon)\right)$ (blue), the out-of-sample performance of variance $V\left(\widehat{x}_{N}^{\mathrm{Var}_{\mathrm{cst}}}(\varepsilon)\right)$ (orange), and the Conditional Value at Risk $\mathrm{CVaR}_{\alpha,\xi\sim\mathbb{P}}\left(-\left\langle \widehat{x}_{N}^{\mathrm{CVaR}_{\mathrm{cst}}}(\varepsilon),\xi\right\rangle\right)$ (blue), as well as Sharpe Ratios as functions of the Wasserstein radius $\varepsilon$. These are estimated based on 200 simulations for $N=300$. In this instance, $\mu=0.15$ and $\alpha=0.05$.
} \label{fig:SimulatedDataCasoConstrainSatisfiesAlways}
\end{figure}

Finally, it is important to mention that the behavior described in Figures \ref{fig:RetursnVsEpsilon} and \ref{fig:RetursnVarianceOptValueVsEpsilon} is influenced by the level of return required in those simulations, where that level was $\mu=0.25$. However, if the minimum required return is lowered, for example, to $\mu=0.15$, then the behavior shown in Figures \ref{fig:RetursnVsEpsilon} and \ref{fig:RetursnVarianceOptValueVsEpsilon} can change considerably, especially in the case where the objective function is CVaR. In this case, the constraint tends to be satisfied with high probability out-of-sample for all values of $\varepsilon$ that make the problem feasible, as illustrated in Figure \ref{fig:SimulatedDataCasoConstrainSatisfiesAlways}. Hence, satisfying the constraint should no longer be the priority, and instead, priority should be given to finding $\varepsilon$ values that minimize the out-of-sample objective function. However, it is important to emphasize that this phenomenon occurs only with the problem that has CVaR as the objective function and only when the required returns are not excessively high. Additionally, this phenomenon is not counter-intuitive, as CVaR and variance are different measures of risk in their behavior. Specifically, CVaR focuses on making the negative tail of $\langle x,\xi\rangle$ as close to zero as possible, which is not necessarily achieved with a small variance. In fact, this could be achieved with high returns, as evidenced in Figure \ref{fig:SimulatedDataCasoConstrainSatisfiesAlways}.

\paragraph{\textbf{Performance of expected confidence level of $\varepsilon$.}}

\begin{figure}[t] 
 \centering
 \begin{tabular}{cc}
    \includegraphics[scale=0.4]{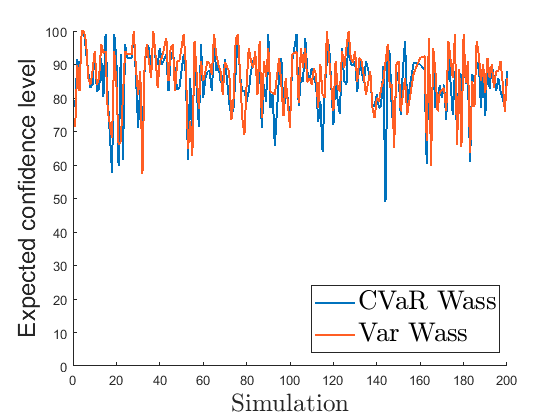} &  \includegraphics[scale=0.4]{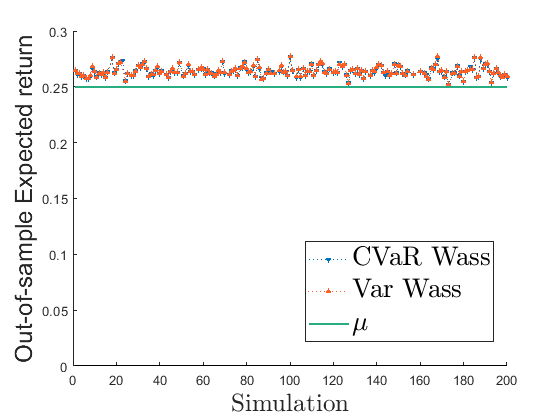}  \\
     (a)  & (b)  
 \end{tabular}
\caption{In (a), the expected confidence level of $\varepsilon=2\widehat{\varepsilon}_{N}^{\mathrm{max}}(\mu)/5$. In (b), the out-of-sample performance of expected returns for $\varepsilon=2\widehat{\varepsilon}_{N}^{\mathrm{max}}(\mu)/5$. This is for each of the problems (\ref{MarkovizDROWReformLargaCVaR}) and (\ref{MarkovizDROWReformDeFormLargaPortflio}), as well as for each of the 200 simulations that were conducted. In this instance, $N=300$ and $\mu=0.25$. } \label{fig:ConfidenceLevel}
\end{figure}

In this part, the aim is to evaluate the performance of the strategy described in subsection \ref{Subsec:Bootstrap}, specifically the expected confidence level of any given $\varepsilon$. To achieve this, we consider $\mu=0.25$ and perform 200 simulations of size $N = 300$. Using $\varepsilon=2\widehat{\varepsilon}_{N}^{\mathrm{max}}(\mu)/5$ and proceeding as in the previous application, Figure \ref{fig:ConfidenceLevel}(a) shows that the expected confidence level is, on average, around 86.7\% and 84.6\% for portfolios generated by solving (\ref{MarkovizDROWReformLargaCVaR}) and (\ref{MarkovizDROWReformDeFormLargaPortflio}), respectively. Furthermore, Figure \ref{fig:ConfidenceLevel}(b) shows that in 100\% of the simulations, the constraint was satisfied for portfolios generated by each problem (\ref{MarkovizDROWReformLargaCVaR}) and (\ref{MarkovizDROWReformDeFormLargaPortflio}), respectively. This suggests that the concept of expected confidence level can be considered a conservative estimator of the probability that the portfolio will satisfy the out-of-sample constraint. 

\subsubsection{Using real market data}

We now turn our attention to real market data. The data utilized in this study comprise the daily returns of 23 companies selected from the S\&P 500 index. The selected returns pertain to the companies described below.

\begin{table}[h]
\begin{center}
\begin{tabular}{lll}
AAPL - Apple & INTC - Intel & PG - P\&G \\
AMZN - Amazon & JNJ - Johnson \& Johnson  & T - AT\&T\\
BAC - Bank of America & JPM - J.P Morgan &  UNH -UnitedHealth Group\\ 
BRKA - Berkshire Hathaway & KO - Coca Cola& VZ - Verizon\\
CVX - Chevron &  MA - Mastercard & WFC - Wells Fargo\\
DIS - Disney & MRK - Merck \& Co & WMT - Walmart\\
HD - The Home Depot & XOM - Exxom Mobil &  MSFT - Microsoft \\ 
PFE - Pfizer  & GOOG - Alphabet Google & 
\end{tabular}
\end{center}
\end{table}

The data encompass the time window from January 1, 2008, to June 30, 2021. In our experiments, we aim to analyze cumulative wealth over time by employing a rolling horizon procedure with daily rebalancing. To illustrate, we use the data from January 1, 2008, to February 13, 2018, to estimate the portfolio vector for February 14, 2018. Subsequently, we use the data from January 2, 2008, to February 14, 2018, to estimate the portfolio vector for February 15, 2018, and so on. This process continues by removing the earliest return and adding a return for the next period until we reach the end of the dataset. Our objective is to observe how cumulative wealth evolves during this time period. Moreover, we compare our approach with standard portfolio optimization techniques, including CVaR SAA, Var SAA, EW, MinCVaR, MinVar, and MaxSR. The techniques labeled SAA involve solving (\ref{StochsticProgWithExpectConsCVaR}) and (\ref{StochsticProgWithExpectConstVariance})  respectively using a Sample Average Approximation approach. The other four techniques are described below:
\begin{enumerate}
    \item[$\bullet$] Equal Weight (EW): This approach assigns equal weight to all assets in the portfolio.
    \item[$\bullet$] Minimum CVaR (MinCVaR): This technique employs the empirical distribution $\widehat{\mathbb{P}}_{N}$ to find the portfolio $x$ that minimizes the expression $\mathrm{CVaR}_{\alpha,\xi\sim\widehat{\mathbb{P}}_{N}}\left(-\langle x,\xi \rangle\right)$ with $\alpha=0.05$.
    \item[$\bullet$] Minimum variance (MinVar): This technique utilizes the sample covariance matrix $\widehat{\Sigma}_{N}$ to find the portfolio $x$ that minimizes the expression $\langle x, \widehat{\Sigma}_{N} x \rangle$.
    \item[$\bullet$] Maximum Sharpe ratio (MaxSR): This technique employs the sample mean vector $\widehat{\mathbf{m}}_{N}$ and the sample covariance matrix $\widehat{\Sigma}_{N}$ to find the portfolio $x$ that maximizes the Sharpe Ratio, defined as $\frac{\langle \widehat{\mathbf{m}}_{N},x \rangle}{\langle x, \widehat{\Sigma}_{N} x \rangle}$.
\end{enumerate}

For a given $\mu$, we focus on analyzing the strategies in (\ref{MarkovizDROWReformLargaCVaR}) and (\ref{MarkovizDROWReformDeFormLargaPortflio}). We will refer to these strategies as CVaR Wass and Var Wass, respectively. Each strategy is evaluated at different epsilon values, specifically $\varepsilon=\widehat{\varepsilon}_{N}^{\mathrm{max}}(\mu)$, $\varepsilon=\frac{3\widehat{\varepsilon}_{N}^{\mathrm{max}}(\mu)}{4}$, and $\varepsilon=\frac{\widehat{\varepsilon}_{N}^{\mathrm{max}}(\mu)}{2}$. To differentiate the strategies concerning the $\varepsilon$ used, we add the terms MaxFact, 3MaxFact/4, and MaxFact/2 to the strategy name, respectively. For instance, if the CVaR Wass strategy is analyzed with $\varepsilon=\widehat{\varepsilon}_{N}^{\mathrm{max}}(\mu)$, it is referred to as CVaR Wass MaxFact. With these conventions established, we present the results of our numerical experiments. Additionally, in the numerical experiments, the minimum allowable daily expected return level $\mu$ was set to 0.001.

\begin{figure}[t] 
 \centering
 \begin{tabular}{c}
     \includegraphics[scale=0.39]{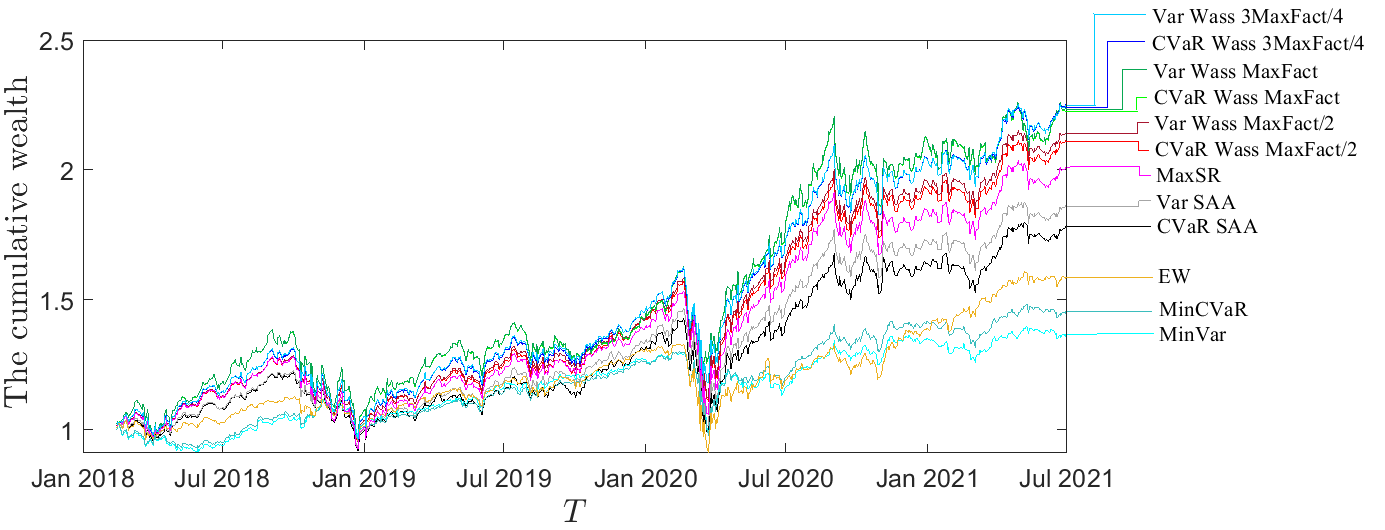}
 \end{tabular}
\caption{The cumulative wealth of the trading strategies with $\mu=0.001$. } \label{fig:CumulativeAndPrtfolioValue} 
\end{figure}

Figure \ref{fig:CumulativeAndPrtfolioValue} displays the cumulative wealth generated by the portfolios resulting from the strategies previously described. The distributionally robust approach based on Wasserstein distances, as proposed in this paper, leads to the highest cumulative wealth, surpassing traditional strategies such as SAA, EW, MinCVaR, MinVar, and MaxSR. Notably, the evaluation period includes the onset of the COVID-19 pandemic, which impacted all investment strategies. Nevertheless, our proposed strategies appear to mitigate the long-term effect on portfolio value.

\begin{table}[tb]
\begin{center}
\begin{tabular}{lllllll} \cline{2-7} 
  & \multicolumn{1}{c}{\multirow{2}{*}{Mean}} &\multicolumn{1}{c}{Standard} & \multicolumn{1}{c}{Sharpe}& \multicolumn{1}{c}{\multirow{2}{*}{Turnover}}  & \multicolumn{1}{c}{Avg. Portfolio}  & \multicolumn{1}{c}{\multirow{2}{*}{$\mathrm{CVaR}_{\alpha=0.05}$}} \\ 
  &  &\multicolumn{1}{c}{deviation} &\multicolumn{1}{c}{Ratio} &   & \multicolumn{1}{c}{Assets} &      \\ \hline
CVaR Wass & \multirow{2}{*}{ 0.0010932} & \multirow{2}{*}{ 0.01719 } & \multirow{2}{*}{0.063594 } & \multirow{2}{*}{ 0.045383 } & \multirow{2}{*}{ 4.3176 } & \multirow{2}{*}{ 0.040446}  \\ 
MaxFact &  & &  &  &  &   \\ \noalign{\global\arrayrulewidth=0.2mm}
  \arrayrulecolor{Silver}\hline
CVaR Wass & \multirow{2}{*}{ 0.0010825} & \multirow{2}{*}{  0.01599 } & \multirow{2}{*}{ 0.067698} & \multirow{2}{*}{ 0.035326 } & \multirow{2}{*}{ 6.2682  } & \multirow{2}{*}{0.037618 }  \\
3MaxFact/4 &  & &  &  &  &   \\ \noalign{\global\arrayrulewidth=0.2mm}
  \arrayrulecolor{Silver}\hline
CVaR Wass & \multirow{2}{*}{ 0.0010047} & \multirow{2}{*}{ 0.015693 } & \multirow{2}{*}{ 0.064023} & \multirow{2}{*}{ 0.040298} & \multirow{2}{*}{ 7.3953} & \multirow{2}{*}{  0.036657}  \\
MaxFact/2 &  & &  &  &  &   \\ \noalign{\global\arrayrulewidth=0.2mm}
  \arrayrulecolor{Silver}\hline
Var Wass & \multirow{2}{*}{0.0010947 } & \multirow{2}{*}{ 0.01719 } & \multirow{2}{*}{ 0.063679} & \multirow{2}{*}{ 0.045387} & \multirow{2}{*}{ 4.2953} & \multirow{2}{*}{ 0.040437}  \\
MaxFact &  & &  &  &  &   \\ \noalign{\global\arrayrulewidth=0.2mm}
  \arrayrulecolor{Silver}\hline
Var Wass & \multirow{2}{*}{ 0.0010793} & \multirow{2}{*}{ 0.015938} & \multirow{2}{*}{ 0.06772} & \multirow{2}{*}{ 0.033119} & \multirow{2}{*}{ 5.9965} & \multirow{2}{*}{ 0.037457}  \\
3MaxFact/4 &  & &  &  &  &   \\ \noalign{\global\arrayrulewidth=0.2mm}
  \arrayrulecolor{Silver}\hline
Var Wass & \multirow{2}{*}{ 0.0010202} & \multirow{2}{*}{ 0.01564} & \multirow{2}{*}{ 0.065232 } & \multirow{2}{*}{ 0.032776 } & \multirow{2}{*}{ 7.2988 } & \multirow{2}{*}{ 0.036351}  \\
MaxFact/2 &  & &  &  &  &   \\ \noalign{\global\arrayrulewidth=0.2mm}
  \arrayrulecolor{Silver}\hline
CVaR SAA &  0.00080191 & 0.0153 & 0.052413 & 0.065448 & 6.8388 & 0.03512    \\ \noalign{\global\arrayrulewidth=0.2mm}
  \arrayrulecolor{Silver}\hline
Var SAA &  0.00085056 &  0.01522 & 0.055883 & 0.040939 & 8.0188 & 0.034432    \\ \noalign{\global\arrayrulewidth=0.2mm}
  \arrayrulecolor{Silver}\hline
MinCVaR & 0.00051007 &  0.011487 & 0.044403  & 0.015437 & 13.711 & 0.026447     \\ \noalign{\global\arrayrulewidth=0.2mm}
  \arrayrulecolor{Silver}\hline
MinVar & 0.00043373  & 0.011263 & 0.03851  & 0.0096731  &  11.331 &  0.026286    \\ \noalign{\global\arrayrulewidth=0.2mm}
  \arrayrulecolor{Silver}\hline
MaxSR & 0.0009563 & 0.016398 & 0.058319 & 0.026923 &  5.3859 & 0.038047   \\ \noalign{\global\arrayrulewidth=0.2mm}
  \arrayrulecolor{Silver}\hline
EW & 0.00063778 & 0.01362 & 0.046828 & 0.97409 & 23 & 0.031383 \\
\hline
\end{tabular}
\end{center}
\caption{Performances of different portfolio strategies.} \label{Tabla:Resultados}
\end{table}

Table \ref{Tabla:Resultados} presents various out-of-sample indicators for different strategies. Note that the mean of all Wasserstein strategies surpasses $\mu$. This observation is significant, particularly since none of the other evaluated strategies, especially the SAA strategies, which serve as the counterpart of the Wasserstein-based strategies, manage to achieve this, despite the large sample size. Indeed, it is generally understood that if the sample size is sufficiently large, the SAA strategies should yield a portfolio close to the one obtained by solving (\ref{StochsticProgWithExpectConst}) if the distribution of returns were known. As a result, this portfolio would satisfy the constraint $\mathbb{E}_{\xi\sim\mathbb{P}}[\left\langle x,\xi\right\rangle]\geq \mu$ with high probability. However, in this case, achieving this outcome may prove more difficult, given that not all data in the sample originate from the same distribution. This situation implies that the mean of returns obtained with the SAA strategies does not exceed $\mu$. In contrast, the Wasserstein-based strategies successfully overcome this challenge.

Another crucial indicator to consider is the standard deviation. In this regard, the standard deviations of the Wasserstein approaches are among the largest. However, the difference relative to the SAA approach is not particularly significant, especially for MaxFact/2 approaches. In addition, the Sharpe Ratios of the Wasserstein-based strategies are the highest, with all of them exceeding 0.063, while the Sharpe Ratios of the other strategies do not surpass 0.059. This observation suggests that the balance between variance and expected value is more favorable for the Wasserstein-based strategies.

Regarding other indicators, it can be observed that turnover increases with $\varepsilon$, though in some instances, it remains lower than that achieved by SAA approaches. Turnover, a metric measuring the percentage of wealth traded when implementing a strategy, is defined in this paper as outlined in \cite{Kang2019}. Generally, lower turnover is preferable due to its impact on the transaction costs associated with the applied strategy. Observations from Table \ref{Tabla:Resultados} reveal that, for CVaR, all Wasserstein-based strategies yield lower turnover compared to their SAA counterparts. Moreover, in the context of variance, only the 3MaxFact/4 and MaxFact/2 Wasserstein-based strategies provide lower turnover than their corresponding SAA versions. Concerning other strategies, all exhibit lower turnover than that generated by Wasserstein-based strategies, with the exception of MaxSR. However, it is crucial to emphasize that these strategies tend to be more conservative compared to Wasserstein-based strategies.

Table \ref{Tabla:Resultados} also presents the average number of assets in the portfolio, referred to as average portfolio assets. This metric reveals the average number of assets constituting the portfolios generated by each strategy during the evaluation period. This information is significant since daily rebalancing in the experiment results in daily changes to the portfolio composition. Furthermore, this indicator can be helpful in identifying the most promising companies within the portfolio. The table shows that the average portfolio assets decrease when $\varepsilon$ increases for Wasserstein-based strategies. For variance-focused Wasserstein-based strategies, the value of this metric is lower than that of the corresponding SAA version. The same holds true for Wasserstein-based strategies centered on CVaR, except for MaxFact/2. When considering the remaining strategies, all exhibit higher average portfolio assets. This observation may be attributed to the more conservative nature of these strategies, which leads them to seek greater diversification within the portfolio.

The final indicator displayed in Table \ref{Tabla:Resultados} is the CVaR at $\alpha=0.05$ generated by each of the strategies under examination. In this regard, Wasserstein-based strategies exhibit the highest CVaR values compared to all other strategies, with the exception of MaxSR. Nonetheless, the difference in value relative to the SAA strategies can be considered not particularly significant, particularly for the MaxFact/2 strategies. This observation can be justified by the fact that Wasserstein-based strategies deliver a higher average return while maintaining a favorable balance with respect to variance. Such behavior may come at a cost, which is manifested in the CVaR.

\begin{figure}[t] 
 \centering
 \begin{tabular}{cc}
    \includegraphics[scale=0.4]{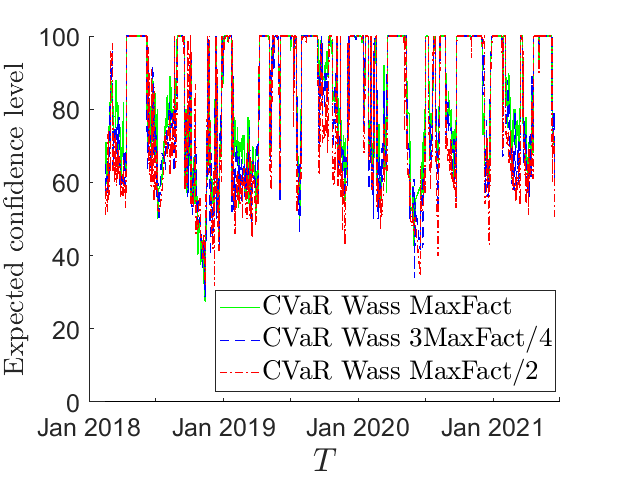} &  \includegraphics[scale=0.4]{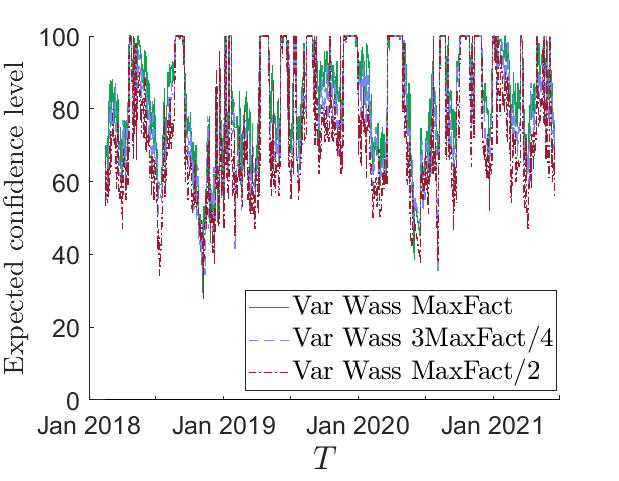}  \\
     (a) Conditional Value at Risk case.  & (b) Variance case.
 \end{tabular}
\caption{ Daily Expected Confidence level for $\mu=0.001$.  } \label{fig:ConfidenceLevelRealData}
\end{figure}

Another aspect warranting discussion is the concept of expected confidence level, introduced in subsection \ref{Subsec:Bootstrap}. In summary, based on the results obtained from simulated data, this concept can be interpreted as a lower estimate of the probability that the constraint will be satisfied out-of-sample. In practice, knowing this information is important for investors since it allows them to anticipate the performance of strategies before making decisions based on its suggestions. In this context, Figure \ref{fig:ConfidenceLevelRealData} displays the daily expected confidence level induced by the Wasserstein strategies. The figure reveals that this level remains above 40\% almost always.

However, when examining the averages, it is noted that for Wasserstein-based strategies focused on CVaR, the average percentages were 84.9\%, 83.3\%, and 81.2\% for MaxFact, 3MaxFac/4, and MaxFact/2, respectively. Similarly, for the variance-focused Wasserstein-based strategies, the average percentages were 83.8\%, 80.3\%, and 75.3\% for MaxFact, 3MaxFac/4, and MaxFact/2, respectively. With this information in mind, it is necessary to determine if these percentages were reflected in reality. This can be assessed by examining Table \ref{Tabla:Resultados}, specifically the column providing information about the means of the returns. As previously mentioned, all Wasserstein-based strategies exceeded the minimum required expected return, which in this case was $\mu=0.001$. Thus, it seems that the concept of expected confidence level does provide insight into the out-of-sample performance of the strategy.

\begin{figure}[t] 
 \centering
 \begin{tabular}{c}
     \includegraphics[scale=0.45]{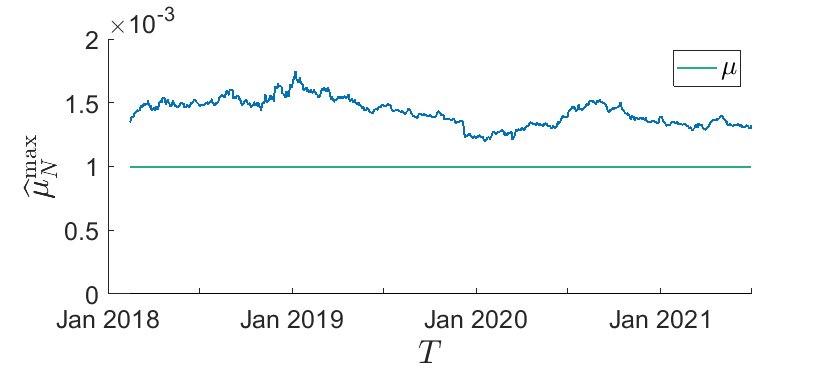}
 \end{tabular}
\caption{Daily maximum expected return level that guarantees feasibility. In this case, $\mu=0.001$. } \label{fig:MuMaxFactibleDatosReales} 
\end{figure}

Finally, a concern that may arise is the possibility of Wasserstein-based strategies admitting a higher minimum expected return than the one utilized in this specific case. Recall that the experiments were conducted with $\mu=0.001$. To address this question, it is crucial to consider the concept of $\widehat{\mu}_{N}^\text{max}$, as defined in Corollary \ref{Corol:MuAndRadioFactibleMeanVarSupportAcot}. In this regard, the value of $\widehat{\mu}_{N}^\text{max}$ determines the maximum value that can be assigned to $\mu$ such that the two Wasserstein-based approaches remain feasible. Figure \ref{fig:MuMaxFactibleDatosReales} presents the daily $\widehat{\mu}_{N}^\text{max}$ values, suggesting that the value of $\mu$ could have been increased.

Nonetheless, it is important to remember that elevating this level also impacts the size of the set of $\varepsilon$ values rendering the Wasserstein-based strategies feasible. This, in turn, affects the expected confidence level. The rationale behind this is that reducing the proximity between $\mu$ and $\widehat{\mu}_{N}^\text{max}$ implies that the $\varepsilon$ values induced by $\mu$, which make Wasserstein-based strategies feasible, possess a lower expected confidence level. Consequently, it becomes more challenging to ensure that the portfolios generated by these strategies satisfy the constraint with a high probability. However, based on the discussions thus far, $\widehat{\mu}_{N}^\text{max}$ and the expected confidence level can be considered tools that enable investors to establish a trade-off between the minimum expected return ($\mu$) and the feasibility of the strategies before making decisions.

\section{Conclusions and future work}  \label{Sec:Conclusions}

In this study, we demonstrated that the Wasserstein distance-based approach (\ref{MarkovizRobust}) possesses an equivalent finite-dimensional formulation for instances where the objective function is either the variance or the Conditional Value at Risk. Moreover, in certain cases, this formulation can be convex. We also established theoretical results characterizing the values of $\mu$ and $\varepsilon$ for which the Wasserstein approach (\ref{MarkovizRobust}) is valid and feasible. For future research, we aim to extend the results presented in this work to a broader set of functions $F$ and $G$, beyond the set determined by Lipschitz functions with respect to $\xi$. Furthermore, we plan to explore alternative types of $\Phi$ functions, such as considering $\Phi$ as a probability. Specifically, we intend to investigate the case where $\Phi(F(x,\xi),\xi\sim\mathbb{P})=\mathbb{P}(F(x,\xi)\geq \rho)$, with $\rho$ being a fixed parameter.

In addition, we applied our strategy to portfolio optimization. In the experiments conducted within this context, it was observed that the proposed strategies (\ref{MarkovizDROWReformLargaCVaR}) and (\ref{MarkovizDROWReformDeFormLargaPortflio}) exhibit performance in line with the priority of satisfying the constraint with a high probability while not sacrificing the out-of-sample value of the objective function too much. This is further supported by the behavior of the Sharpe Ratios observed in the experiments using real data and the high levels of accumulated wealth compared to those generated by other evaluated strategies. Moreover, it became clear that the concept of expected confidence level proposed in this paper can be viewed as a lower estimate of the probability that the portfolios generated by the proposed strategies satisfy the constraint of (\ref{StochsticProgWithExpectConsCVaR}) and (\ref{StochsticProgWithExpectConstVariance}) out-of-sample. This positions this concept as a valuable tool for forecasting the performance of portfolios created by the proposed strategies before making a decision and for calibrating the size of the ambiguity set.

\textbf{Acknowledgments}\\
This work was supported by the Research Fund of the Faculty of Sciences of the Universidad de los Andes INV-2021-128-2307 and INV-2021-126-2273.

\appendix

\section{Proofs of Lemmas and Theorems}

We present proofs of the results presented in this work. Section \ref{Apendice:PruebasLema} present the proof of Lemma \ref{Lemma:Confianza}, and  Section \ref{Sec:Appendix:ProofThmMeanVarianze} shows the proofs of Theorems \ref{Prop:Reformul1MarkovizDROWReformLargaLargaSupportAcot} and \ref{Prop:Reformul1MarkovizDROWReformLargaLarga} and its corollaries.

\subsection{Proofs of Lemma \ref{Lemma:Confianza}}\label{Apendice:PruebasLema}

For the proof of Lemma \ref{Lemma:Confianza}, it is necessary the following result which was proved in \cite{fonsecaDecDep2023}.

\begin{lemma} \label{Lemma:JustificacionBola}
Assuming that $F$ satisfies Assumption \ref{AssumptionPrincipal}, $ W_{p}(\widehat{\mathbb{P}}_{N}^{x,F},\mathbb{P}^{x,F}) \leq \gamma_{x,F,q}  W_{p}\left( \widehat{\mathbb{P}}_{N} ,\mathbb{P}\right)$ for $p,q\geq1$ where $ W_{p}\left( \widehat{\mathbb{P}}_{N} ,\mathbb{P}\right)$ is considered with cost function $\mathbf{d}=\|\cdot\|_{q}$ and $W_{p}(\widehat{\mathbb{P}}_{N}^{x,F},\mathbb{P}^{x,F}) $ with cost function $\mathbf{d}=|\cdot|$.
\end{lemma}

\proof[Lemma \ref{Lemma:Confianza}]
Let $\varepsilon>0$ such that $\mathbb{P}\in\mathcal{B}_{\varepsilon}(\widehat{\mathbb{P}}_{N})$ and (\ref{MarkovizRobust}) is feasible, then, to simplify the notation in this proof, we denote $\widehat{x}_{N,p,q}^{A,cst}(\varepsilon)$ by $\widehat{x}_{N}(\varepsilon)$. Therefore, by Lemma \ref{Lemma:JustificacionBola}, we have that $\mathbb{P}^{\widehat{x}_{N}(\varepsilon),G}\in\mathcal{B}_{\varepsilon\gamma_{\widehat{x}_{N}(\varepsilon),G,q} }\left(\widehat{\mathbb{P}}_{N}^{\widehat{x}_{N}(\varepsilon),G}\right)$. Hence, because $\inf_{\mathbb{Q}\in\mathcal{B}_{\varepsilon\gamma_{\widehat{x}_{N}(\varepsilon),G,q}}\left(\widehat{\mathbb{P}}_{N}^{\widehat{x}_{N}(\varepsilon),G}\right)}\mathbb{E}_{\xi\sim\mathbb{Q}}[\zeta]\geq \mu$, we obtain $\mathbb{E}_{\xi\sim\mathbb{P}^{\widehat{x}_{N}(\varepsilon),G}}[\zeta]\geq \mu$. However, note that $\mathbb{E}_{\xi\sim\mathbb{P}^{\widehat{x}_{N}(\varepsilon),G}}[\zeta]=\mathbb{E}_{\xi\sim\mathbb{P}^{\widehat{x}_{N}(\varepsilon),G}}\left[\zeta^{\widehat{x}_{N}(\varepsilon),G}\right]=\mathbb{E}_{\xi\sim\mathbb{P}}\left[G\left( \widehat{x}_{N}(\varepsilon),\xi\right)\right]$. Therefore, we conclude that  $\mathbb{E}_{\xi\sim\mathbb{P}}\left[G\left( \widehat{x}_{N}(\varepsilon),\xi\right)\right] \geq \mu$. 
\qed
\endproof

\subsection{Proofs of Theorems \ref{Prop:Reformul1MarkovizDROWReformLargaLargaSupportAcot} and \ref{Prop:Reformul1MarkovizDROWReformLargaLarga}, and Corollary \ref{Corol:MuAndRadioFactibleMeanVarSupportAcot}} \label{Sec:Appendix:ProofThmMeanVarianze}

\proof[Theorem \ref{Prop:Reformul1MarkovizDROWReformLargaLargaSupportAcot}]
Let $\mathbb{X}$  the feasible set of (\ref{MarkovizRobust}), then, by Theorem 2-(a) in \cite{fonsecaDecDep2023}, we have that
\begin{equation*}
\mathbb{X}=\left\{x\in\mathbb{R}^{m}\:\left|\: \max\left\{{\displaystyle \frac{1}{N}\sum_{i=1}^{N}G\left( x,\widehat{\xi}_{i}\right)-\varepsilon\gamma_{x,G}},A_{G}(x) \right\}  \geq \mu,\:  x\in\mathcal{X}\:\right. \right\}.  \label{Set:XcaracterizacionSupporAcot}
\end{equation*}
Therefore, (\ref{MarkovizRobust}) is equivalent to 
\begin{equation*}  
   \underset{ x\in\mathbb{X}  }{\min}  \sup_{\mathbb{Q}\in\mathcal{B}_{\varepsilon \gamma_{x,F}}(\widehat{\mathbb{P}}_{N}^{x,F}) }  \mathbb{E}_{\mathbb{Q}}\left[\zeta\right],
\end{equation*}
but, again,  by Theorem 2-(a) in \cite{fonsecaDecDep2023}, we have that
\begin{align*}
    \underset{ x\in\mathbb{X}  }{\min} \min\left\{\frac{1}{N}{\displaystyle\sum_{i=1}^{N}}F\left(x,\widehat{\xi}_{i}\right)+\varepsilon\gamma_{x,F},B_{F}(x)\right\},
\end{align*}
which is equivalent to (\ref{MarkovizDROWReformLargaSupportAcot}). Analogously, to prove (\ref{MarkovizDROWReformLargaSupportNoAcot}), we use  Theorem 2-(b)  in \cite{fonsecaDecDep2023} Lemma.  
\qed
\endproof

\proof[Theorem \ref{Prop:Reformul1MarkovizDROWReformLargaLarga}]
Let $\mathbb{X}$  the feasible set of (\ref{MarkovizRobust}), then, by Theorem 2-(a) in \cite{fonsecaDecDep2023}, we have that
\begin{align*}
\mathbb{X}= \left\{x\in\mathbb{R}^{m}\:\left|\: \frac{1}{N}\sum_{i=1}^{N}G\left( x,\widehat{\xi}_{i}\right)-\varepsilon\gamma_{x,G} \geq \mu,\:  x\in\mathcal{X}\:\right. \right\}. 
\end{align*}
Therefore, (\ref{MarkovizRobust}) is equivalent to 
\begin{align}
\underset{ x\in\mathbb{X}  }{\min}  \sup_{\mathbb{Q}\in\mathcal{B}_{\varepsilon \gamma_{x,F}}(\widehat{\mathbb{P}}_{N}^{x}) }  \mathrm{Var}_{\mathbb{Q}}\left[\zeta\right]  \label{MarkovizDROWReformLarga1Previa}  
\end{align}  
Additionally, by Theorem 3 in \cite{fonsecaDecDep2023}, (\ref{MarkovizDROWReformLarga1Previa} ) can be rewritten as
\begin{align}
\widehat{J}_{N} (\varepsilon) &=\underset{ x\in\mathbb{X}  }{\mathrm{minimize}} \left( \sqrt{\frac{1}{N}{\displaystyle\sum_{i=1}^{N}}F\left(x,\widehat{\xi}_{i}\right) ^{2}-\frac{1}{N^{2}}\left({\displaystyle\sum_{i=1}^{N}}F\left( x,\widehat{\xi}_{i}\right)\right)^{2} }+\varepsilon\gamma_{x,F} \right)^{2} \nonumber \\
&= \left\{\begin{array}{ll}\underset{ x\in\mathbb{R}^{m}  }{\mathrm{minimize}} &  \left( \sqrt{\frac{1}{N}{\displaystyle\sum_{i=1}^{N}}F\left(x,\widehat{\xi}_{i}\right) ^{2}-\frac{1}{N^{2}}\left({\displaystyle\sum_{i=1}^{N}}F\left( x,\widehat{\xi}_{i}\right)\right)^{2} }+\varepsilon\gamma_{x,F} \right)^{2} \\[0.4cm] \mbox{subject to} & {\displaystyle \frac{1}{N}\sum_{i=1}^{N}G\left( x,\widehat{\xi}_{i}\right)-\varepsilon\gamma_{x,G} \geq \mu,} \\[0.4cm] 
&  x\in\mathcal{X}. \end{array}\right. \nonumber
\end{align} 
\qed
\endproof

\proof[Corollary  \ref{Corol:MuAndRadioFactibleMeanVarSupportAcot}]
 To determine feasibility, it suffices to find an $x\in\mathcal{X}$ such that the constraint is satisfied. We begin with case (i). The hypotheses for this case indicate the existence of an $x\in\mathcal{X}$ such that  $\inf_{\xi\in\Xi}G(x,\xi)\geq \mu$. This implies that $\max\left\{\frac{1}{N}\sum\limits_{i=1}^{N}G(x,\widehat{\xi}_{i})-\varepsilon\gamma_{x,G,q},\inf_{\xi\in\Xi}G(x,\xi)\right\} \geq \mu$ for all $\varepsilon>0$, concluding the proof for this case.

For case (ii), the hypotheses state that ${\displaystyle\sup_{x\in\mathcal{X}}\inf_{\xi\in\Xi} }G(x,\xi)< \mu$, which implies the necessity of finding an $x\in\mathcal{X}$ satisfying $\frac{1}{N}\sum\limits_{i=1}^{N}G(x,\widehat{\xi}_{i})-\varepsilon\gamma_{x,G,q}\geq \mu$. Given that  $\mu < \widehat{\mu}_{N}^{\mathrm{max}}$, there exists $x\in\mathcal{X}$ such that $\frac{1}{N}\sum\limits_{i=1}^{N}G(x,\widehat{\xi}_{i})>\mu$, implying $\frac{\frac{1}{N}\sum\limits_{i=1}^{N}G(x,\widehat{\xi}_{i})-\mu}{\gamma_{x,G,q}}>0$. This further implies $\widehat{\varepsilon}_{N,p,q}^{\mathrm{max}}(\mu)>0$. Therefore, since $\varepsilon<\widehat{\varepsilon}_{N,p,q}^{\mathrm{max}}(\mu)$, there exists $x'\in\mathcal{X}$ such that $\frac{\frac{1}{N}\sum\limits_{i=1}^{N}G(x',\widehat{\xi}_{i})-\mu}{\gamma_{x',G,q}}\geq \varepsilon$, allowing us to conclude $\frac{1}{N}\sum\limits_{i=1}^{N}G(x',\widehat{\xi}_{i})- \varepsilon \gamma_{x',G,q} \geq \mu$.

Lastly, for case (iii), the proof follows a similar structure to case (ii) but with $p\neq 1$.
\qed
\endproof

\bibliographystyle{splncs04}
\bibliography{refs}
\end{document}